\documentclass{cmslatex}
\usepackage[utf8]{inputenc}
\usepackage[paperwidth=7in, paperheight=10in, margin=.875in]{geometry}
\usepackage[backref,colorlinks,linkcolor=red,anchorcolor=green,citecolor=blue]{hyperref}
\usepackage{booktabs,multirow}
\usepackage[noend]{algpseudocode}
\usepackage[pagewise]{lineno}
\usepackage{algorithmicx,algorithm,setspace,subfigure}
\usepackage{amsfonts,amssymb}
\usepackage{amsmath}
\usepackage{graphicx}
\usepackage{amsmath}
\usepackage{cite}
\usepackage{subfigure}
\usepackage{todonotes}
\usepackage{enumerate}
\usepackage{ulem,marvosym}

\usepackage[T1]{fontenc} 

\renewcommand{\theequation}{\arabic{section}.\arabic{equation}}

\newcommand{\abs}[1]{\left|#1\right|}

\newcommand{\mbr}{\mathbb{R}}
\newcommand{\bu}{\boldsymbol{u}}
\newcommand{\bv}{\boldsymbol{v}}
\newcommand{\bgamma}{\boldsymbol{\gamma}}
\newcommand{\balpha}{\boldsymbol{\alpha}}
\newcommand{\bbeta}{\boldsymbol{\beta}}
\newcommand{\bphi}{\boldsymbol{\phi}}
\newcommand{\bpsi}{\boldsymbol{\psi}}

\newcommand{\bw}{\boldsymbol{w}}
\newcommand{\bC}{\boldsymbol{C}}

\newcommand{\bJ}{\boldsymbol{J}}
\newcommand{\bP}{\boldsymbol{P}}
\newcommand{\bQ}{\boldsymbol{Q}}

\newcommand{\bvarphi}{\boldsymbol{\varphi}}

\newcommand{\bbR}{\mathbb{R}}

\floatname{algorithm}{Algorithm}

\let\oldequation\equation
\let\oldendequation\endequation

\renewenvironment{equation}{\linenomathNonumbers\oldequation}{\oldendequation\endlinenomath}

\let\oldalign\align
\let\oldendalign\endalign

\renewenvironment{align}{\linenomathNonumbers\oldalign}{\oldendalign\endlinenomath}

\allowdisplaybreaks
\begin{document}

\let\oldequationa\equation*
\let\oldendequationa\endequation*

\renewenvironment{equation*}{\linenomathNonumbers\oldequationa}{\oldendequationa\endlinenomath}

\begin{sloppypar}

 \title{
 A numerical algorithm with linear complexity for Multi-marginal Optimal Transport with $L^1$ Cost
 \thanks{Received date, and accepted date (The correct dates will be entered by the editor).}
 }

          \author{Chunhui Chen\thanks{Department of Mathematical Sciences, Tsinghua University, Beijing 100084, China (\href{mailto:cch21@mails.tsinghua.edu.cn}{cch21@mails.tsinghua.edu.cn}).}
          \and Jing Chen\thanks{School of Physical and Mathematical Sciences, Nanyang Technological University, Singapore 639798(\href{mailto:jing.chen@ntu.edu.sg}{jing.chen@ntu.edu.sg}).}
          \and Baojia Luo\thanks{Theory Lab, Central Research Institute, 2012 Labs, Huawei Technologies Co. Ltd., Hong Kong, China(\href{mailto:luobaojia2@huawei.com}{luobaojia2@huawei.com}).}
          \and Shi Jin\thanks{School of Mathematical Sciences, Institute of Natural Sciences, and MOE-LSC Shanghai Jiao Tong University, Shanghai 200240, China (\href{mailto:shijin-m@sjtu.edu.cn}{shijin-m@sjtu.edu.cn}).}
          \and Hao Wu\thanks{Corresponding author. Department of Mathematical Sciences, Tsinghua University, Beijing 100084, China (\href{mailto:hwu@tsinghua.edu.cn}{hwu@tsinghua.edu.cn})}.}

 \pagestyle{myheadings} \markboth{
 }{Authors} \maketitle
\begin{abstract}
    Numerically solving multi-marginal optimal transport (MMOT) problems is computationally prohibitive, even for moderate-scale instances involving $l\ge4$ marginals with support sizes of $N\ge1000$. 
    The cost in MMOT is represented as a tensor with $N^l$ elements. Even accessing each element once incurs a significant computational burden. In fact, many algorithms require direct computation of tensor-vector products, leading to a computational complexity of $O(N^l)$ or beyond. In this paper, inspired by our previous work [\textit{Comm. Math. Sci.,} 20 (2022), pp. 2053 – 2057], we observe that the costly tensor-vector products in the Sinkhorn Algorithm can be computed with a recursive process by separating summations and dynamic programming. Based on this idea, we propose a fast tensor-vector product algorithm to solve the MMOT problem with $L^1$ cost, achieving a miraculous reduction in the computational cost of the entropy regularized solution to $O(N)$. Numerical experiment results confirm such high performance of this novel method which can be several orders of magnitude faster than the original Sinkhorn algorithm.
\end{abstract}



\section{Introduction}

Multi-marginal optimal transport (MMOT), first proposed by Gangbo and \'Swiech \cite{gangbo1998optimal}, is an extension of the classical optimal transport problem \cite{peyre2019computational,benamou2016numerical}. It aims to find an optimal transport plan that minimizes total cost while fitting multiple marginal distributions. 
MMOT problems naturally arise in various fields, such as machine learning \cite{agueh2011barycenters,cao2019multi,metelli2019propagating,haasler2021multi}, incompressible fluid dynamics \cite{benamou2019generalized,baradat2020small}, density function theory \cite{buttazzo2012optimal,cotar2013density,hu2021global}, Schrödinger bridge \cite{DBLP:journals/jota/ChenGP16, DBLP:journals/siamco/HaaslerRCK21}, and tomographic reconstruction \cite{abraham2017tomographic}, and thus have attracted wide attention in recent years.


However, the heavy computational burden of solving general MMOT problems limits its broad application. For $l$-marginal distributions with support sizes of $N$, the $l-$th order cost tensor in MMOT problems contains $N^l$ elements. To fully obtain the information of the cost tensor, it is inevitable to repeatedly access all elements in the tensor, leading to a significant computational cost. For example, the generalized Sinkhorn algorithms \cite{benamou2015iterative,benamou2016numerical,nenna2016numerical,lin2022complexity,tupitsa2020multimarginal} require repeated computation of tensor-vector products, resulting in a computational complexity of $O(N^l)$. More severely, directly solving linear programming problems has a computational complexity of $N^{O(l)}$.
 Therefore, these methods remain computationally prohibitive even for moderate-scale MMOT problems. Some modified algorithms with lower computational complexity have been proposed for specific MMOT problems, such as the MMOT problem with a tree structure \cite{DBLP:journals/siamco/HaaslerRCK21} and the Wasserstein barycenter \cite{benamou2015iterative,altschuler2021wasserstein}.

In this work, we propose a novel implementation of the Sinkhorn algorithm for solving the entropy regularized MMOT problem with $L_1$ cost applications in image processing \cite{solomon2014earth}, computer vision \cite{pele2009fast} and seismic tomography \cite{metivier2016measuring}, which has linear computational complexity relative to support size $N$. 
This work is a follow-up work of the fast Sinkhorn algorithms \cite{liao2022fast,liao2023fast2} , which observe the special structure of the kernel matrix with Wasserstein-1 metric and utilize dynamic programming techniques \cite{horner1819xxi,kleinberg2006algorithm,li2018parallel}, achieving linear computational complexity for solving up to 2-marginal optimal transport problem. Unlike the previous situation, the kernel matrix evolves into an $l$-th order tensor in the $l$-marginal optimal transport problem. The computational burden of the Sinkhorn algorithm becomes prohibitive due to the $O(N^l)$ operations required by the tensor-vector products. 
To address this problem, we observe a similar special structure of the kernel tensor and accelerate the tensor-vector products using the series rearrangement and dynamic programming techniques \cite{kleinberg2006algorithm}, which results in a fast Sinkhorn algorithm with $O(N)$ computational complexity.

The rest of the paper is organized as follows. In Section \ref{sec:mmot}, we review the MMOT problem and the generalized Sinkhorn algorithm for the 1-dimensional (1D) and 3-marginal case. Then, we introduce the key fast tensor-vector product technique and provide a detailed implementation of our algorithm in Section \ref{sec:FS-3}. This algorithm can be conducted in more general scenarios, such as high-dimensional and $l$-marginal optimal transport problems, which are presented in Section \ref{sec:extension}. In Section \ref{sec:numerical experiment}, numerical experiments are carried out to demonstrate the overwhelming advantage of our algorithm in terms of computational efficiency. Finally, we conclude this paper in Section \ref{sec:conclusion}.

\section{Multi-marginal optimal transport and the Sinkhorn algorithm}\label{sec:mmot}

We first review the MMOT problem and the Sinkhorn algorithm \cite{benamou2015iterative}. To streamline our discussion, we exclusively showcase the 3-marginal optimal transport problem and its corresponding Sinkhorn algorithm in 1D space. However, it is important to note that this algorithm can be seamlessly extrapolated to address more marginal cases in higher dimensional space. For a more comprehensive understanding, we refer the readers to \cite{benamou2015iterative}.

We consider the discrete MMOT problem, which is a discretization of the continuous MMOT problem \cite{nenna2016numerical}. It solves the following minimization problem
\begin{equation}\label{eqn:mot problem}
	W(\bu,\bv,\bw)= \inf_{\mathcal{T} \in \Pi}\;\left \langle \mathcal{C}, \mathcal{T} \right \rangle =\inf_{\mathcal{T} \in \Pi}\;\sum\limits_{i=1}^{N_1} \sum\limits_{j=1}^{N_2} \sum\limits_{k=1}^{N_3} c_{ijk} t_{ijk},
\end{equation}
in which $\bu \in \mbr^{N_1},\bv \in \mbr^{N_2},\bw \in \mbr^{N_3}$ are three discrete probabilistic distributions 
\begin{equation*}
	\bu=(u_1,u_2,\cdots,u_{N_1}), \quad \bv=(v_1,v_2,\cdots,v_{N_2}), \quad \bw=(w_1,w_2,\cdots,w_{N_3}).
\end{equation*}
Here $\mathcal{C}=(c_{ijk}) \in \bbR^{N_1\times N_2 \times N_3}$ is the cost tensor, and $\mathcal{T}=(t_{ijk}) \in \bbR_{0+}^{N_1\times N_2 \times N_3}$ is the multi-marginal transport plan, satisfying the linear constraints
\begin{equation}\label{eqn: transport map set}
	\mathcal{T} \in \Pi=\Big\{ \mathcal{T}=(t_{ijk}) \Big| 
	\sum\limits_{j=1}^{N_2} \sum\limits_{k=1}^{N_3} t_{ijk}=u_i,\; 
	\sum\limits_{i=1}^{N_1} \sum\limits_{k=1}^{N_3} t_{ijk}=v_j,\; 
	\sum\limits_{i=1}^{N_1} \sum\limits_{j=1}^{N_2} t_{ijk}=w_k,\;t_{ijk} \ge 0 
	 \Big\}.
\end{equation}
In this work, our discussion is suitable for any $N_1$, $N_2$ and $N_3$. For the sake of simplicity, we assume $N_1 = N_2 = N_3 = N$ in the rest of the paper.

The generalized Sinkhorn algorithm \cite{benamou2015iterative} was proposed to solve the regularized MMOT problem by introducing an entropy regularization term:
\begin{equation*}
	W_\varepsilon(\bu,\bv,\bw)= \inf_{\mathcal{T} \in \Pi}\;\sum\limits_{i,j,k=1}^{N} c_{ijk} t_{ijk} + \varepsilon t_{ijk} \ln{(t_{ijk})}.
\end{equation*}
The above optimization problem can be solved using the method of Lagrange multipliers, with the Lagrangian given by
\begin{multline} \label{eqn:lagrange}
	\mathcal{L}(\mathcal{T},\balpha,\bbeta,\bgamma) =  \sum_{i,j,k=1}^{N} \varepsilon\ln \left(\frac{t_{ijk}}{K_{ijk}} \right)t_{ijk}
	+ \sum_{i=1}^N\alpha_i\left(\sum_{j,k=1}^{N} t_{ijk}-u_i\right) \\
	+ \sum_{j=1}^N\beta_j \left(\sum_{i,k=1}^{N} t_{ijk}-v_j\right)
	+ \sum_{k=1}^N\gamma_k \left(\sum_{i,j=1}^{N} t_{ijk}-w_k\right),
\end{multline}
where $K_{ijk} = e^{\frac{-c_{ijk}}{\varepsilon}}$. Define $\bphi=(\phi_1,\cdots,\phi_N)$, $\bpsi=(\psi_1,\cdots,\psi_N)$, $\bvarphi=(\varphi_1,\cdots,\varphi_N)$ as
\begin{equation*}
	\phi_i=e^{-\frac{1}{3}-\frac{\alpha_i}{\varepsilon}}, \quad \psi_j=e^{-\frac{1}{3}-\frac{\beta_j}{\varepsilon}}, \quad \varphi_k=e^{-\frac{1}{3}-\frac{\gamma_k}{\varepsilon}}.
\end{equation*}
Taking the derivative of the Lagrangian with respect to $t_{ijk}$ and setting it equal to zero yields
\begin{equation}\label{eqn:T value}
	t_{ijk} =  \phi_i \psi_j \varphi_k K_{ijk},\quad  \;1\le i,j,k \le N.
\end{equation}
Combining Equation \eqref{eqn:T value} and the constraints in Equation \eqref{eqn: transport map set}, we can obtain that when the regularized MMOT problem attains its optimal solution, $\bphi, \bpsi, \bvarphi$ satisfy:
\begin{equation}\label{eqn: 3-marginal optimality}
	\phi_i \sum_{j,k=1}^N K_{ijk} \psi_j \varphi_k = u_i, \quad  
	\psi_j \sum_{i,k=1}^N K_{ijk} \phi_i \varphi_k = v_j, \quad
	\varphi_k \sum_{i,j=1}^N K_{ijk} \phi_i \psi_j = w_k. 
\end{equation}
Since the elements in $\mathcal{K}=(K_{ijk})$ are strictly positive, the generalized Sinkhorn algorithm \cite{benamou2015iterative} can be applied to find scaling variables $\bphi,\bpsi,\bvarphi$ satisfying the above equations by iteratively updating
\begin{align}
	\bphi^{(t+1)} &=  \bu \oslash (\mathcal{K} \times_j \bpsi^{(t)} \times_k \bvarphi^{(t)}), \label{eqn:update formula 1} \\
	\bpsi^{(t+1)} &=  \bv \oslash (\mathcal{K} \times_i \bphi^{(t+1)} \times_k \bvarphi^{(t)}), \label{eqn:update formula 2} \\
	\bvarphi^{(t+1)} &=  \bw \oslash (\mathcal{K} \times_i \bphi^{(t+1)} \times_j \bpsi^{(t+1)}), \label{eqn:update formula 3}
\end{align}
in which $\oslash$ represents pointwise division and $t$ denotes the iterative step. 
Here the tensor-vector products follow the notation in \cite{de2000multilinear}, having the forms of
\begin{align}
    & (\mathcal{K} \times_j \bpsi^{(t)} \times_k \bvarphi^{(t)})_i = 
    \sum_{j=1}^N \sum_{k=1}^N K_{ijk} \psi_j^{(t)} \varphi_k^{(t)}, \label{eqn:tensor vector product 1} \\
    & (\mathcal{K} \times_i \bphi^{(t+1)} \times_k \bvarphi^{(t)})_j = 
    \sum_{i=1}^N \sum_{k=1}^N K_{ijk} \phi_i^{(t+1)} \varphi_k^{(t)}, \label{eqn:tensor vector product 2} \\
    & (\mathcal{K} \times_i \bphi^{(t+1)} \times_j \bpsi^{(t+1)})_k = \sum_{i=1}^N \sum_{j=1}^N K_{ijk} \phi_i^{(t+1)} \psi_j^{(t+1)}. \label{eqn:tensor vector product 3}
\end{align}
The final MMOT distance is given by
\begin{equation}\label{eqn:mmot distance}
    W_\epsilon(\bu,\bv,\bw) = \sum_{i,j,k=1}^N \phi_i \psi_j \varphi_k c_{ijk}K_{ijk} = \Big\langle \bvarphi, (\mathcal{C} \odot \mathcal{K}) \times_i \bphi \times_j \bpsi \Big \rangle,
\end{equation}
in which $\odot$ represents the Hadamard product of tensors. The generalized Sinkhorn algorithm is presented in Algorithm \ref{alg:Sinkhorn}.

\medskip

\begin{remark}
Theoretically, obtaining approximate solutions with high accuracy requires an adequately small entropic parameter $\varepsilon$, which may cause over- or underflow in numerical calculations during the Sinkhorn iteration. The log-domain stabilization technique \cite{chizat2018scaling} can be employed to address this issue. When the infinite norm of $\phi$, $\psi$ or $\phi$ exceeds a given threshold $\tau$, the excessive parts of $\phi$, $\psi$, $\phi$ are absorbed into $\balpha$, $\bbeta$, $\bgamma$ to avoid the over- and underflow:
\begin{equation*}
    \begin{array}{lll}
     \balpha \leftarrow \balpha + \varepsilon\ln(\bphi), \quad 
    & \bbeta \leftarrow \bbeta + \varepsilon\ln(\bpsi), \quad
    & \bgamma \leftarrow \bgamma + \varepsilon\ln(\bvarphi), \\
    \bphi \leftarrow \boldsymbol{1}_N, \quad 
    & \bpsi \leftarrow \boldsymbol{1}_N, \quad 
    & \bvarphi \leftarrow \boldsymbol{1}_N.
    \end{array}
\end{equation*}
Correspondingly, the tensor $\mathcal{K}$ should also be rescaled as $K_{ijk} \leftarrow e^{\frac{\alpha_i+\beta_j+\gamma_k}{\varepsilon}}K_{ijk}$.
\end{remark}

\begin{algorithm}  
	\caption{Generalized Sinkhorn Algorithm}  
	\label{alg:Sinkhorn}
	\hspace*{0.02in} {\bf Input:} $\bu, \bv, \bw$ of size $(N,1)$, $\varepsilon$, tol, itr\_max \\
	\hspace*{0.02in} {\bf Output:} $W_\epsilon(\bu,\bv,\bw)$
	\begin{algorithmic}[1] 
		\State initialize 
        $\bphi,\bpsi,\bvarphi = \frac{1}{N}\boldsymbol{1}_{N}$, 
        $t=0$, 
        $\mathrm{Res} = \mathrm{Inf}$
		\While{($t<$ itr\_max) \& ($\mathrm{Res} >$ tol)}
        \State $t \leftarrow t+1$
		\State $\bphi \leftarrow \bu \oslash \big( \mathcal{K} \times_j \bpsi \times_k \bvarphi \big)$
		\State $\bpsi \leftarrow \bv \oslash \big( \mathcal{K} \times_i \bphi \times_k \bvarphi \big)$
		\State $\bvarphi \leftarrow \bw \oslash \big( \mathcal{K} \times_i \bphi \times_j \bpsi \big)$
		\State $\mathrm{Res} \leftarrow \mathrm{sum}\Big (\mid \bphi \odot  \big( \mathcal{K} \times_j \bpsi \times_k \bvarphi \big) -\bu \mid +\mid \bpsi \odot \big( \mathcal{K} \times_i \bphi \times_k \bvarphi \big)-\bv \mid \Big )$
		\EndWhile
          \State \Return{$W_\epsilon(\bu,\bv,\bw)=  \Big\langle \bvarphi, (\bC \odot \mathcal{K}) \times_i \bphi \times_j \bpsi \Big \rangle $}
	\end{algorithmic}  
\end{algorithm}

\section{{An algorithm for MMOT with linear complexity}}\label{sec:FS-3}

The computational complexity of Algorithm \ref{alg:Sinkhorn} is $O(N^3)$ for $3$-marginal MMOT problems, and it scales to $O(N^l)$ for the $l$-marginal case. 
The bottleneck arises from the repeated tensor-vector products related to $\mathcal{K}$ and $\mathcal{C} \odot \mathcal{K}$ in lines 4-8. 

Inspired by the fast Sinkhorn algorithms for Wasserstein-1 distance \cite{liao2022fast,liao2023fast2}, we aim to reduce the computational complexity from $O(N^l)$ to $O(N)$. 
For the sake of simplicity, we discuss the acceleration of Algorithm \ref{alg:Sinkhorn} for $3$-marginal case in this section. The implementation can be naturally extended to the $l$-marginal case for any $l\ge 3$, which will be discussed in Subsection \ref{sec:k-marginal}. On a uniform mesh with a grid spacing of $h$\footnote{The discussion can be naturally extended to the case of non-uniform mesh.}, the cost tensor $\mathcal{C}=(c_{ijk})$ based on the $L^1$ norm writes
\begin{equation*}
	c_{ijk} = (\abs{i-j}+\abs{i-k}+\abs{j-k})h.
\end{equation*}
Thus, the corresponding kernel tensor $\mathcal{K} = (K_{ijk})$ satisfies
\begin{equation*}
	K_{ijk} = \lambda^{\abs{i-j}+\abs{i-k}+\abs{j-k}},
\end{equation*}
in which $\lambda = e^{-\frac{h}{\varepsilon}}$. 
 
To eliminate the absolute value operator in $\mathcal{K}$, we categorize the subscripts $(i,j,k)$ according to their orders into 6 sets\footnote{For the $l$-marginal case, we need to categorize subscripts $(i_1,i_2,\cdots,i_l)$. Therefore, $l!$ sets are required.} 
\begin{align}\label{eqn:D_p}
    \begin{array}{l@{\quad}l}
    D_{1}=\Big\{(i,j,k)\;\Big|\;1\le i \le j \le k \le N \Big\}, &  
	D_{2}=\Big\{(i,j,k)\;\Big|\;1\le j  <  i \le k \le N \Big\},  \\ 
	D_{3}=\Big\{(i,j,k)\;\Big|\;1\le i \le k  <  j \le N \Big\}, & 
	D_{4}=\Big\{(i,j,k)\;\Big|\;1\le j \le k  <  i \le N \Big\},  \\
	D_{5}=\Big\{(i,j,k)\;\Big|\;1\le k  <  i \le j \le N \Big\}, & 
	D_{6}=\Big\{(i,j,k)\;\Big|\;1\le k  <  j  <  i \le N \Big\}. 
   \end{array}
\end{align}
Then, the element of tensor $\mathcal{K}$ has a simplified form of
\begin{equation}\label{eqn:form of K}
    K_{ijk} = \lambda^{a_p i + b_p j + c_p k}, \quad \forall (i,j,k) \in D_{p}, \quad p=1,\cdots,6,
\end{equation}
where $a_p, b_p, c_p$ have specific values
\begin{align}\label{eqn:a_p b_p c_p}
    \begin{array}{ll@{\quad}ll@{\quad}ll@{\quad}ll@{\quad}ll@{\quad}ll}
        a_1 =-&2, \; &a_2 = &0, \; &a_3 =-&2, \; &a_4 = &2, \; &a_5 = &0, \; &a_6 = &2,  \\
        b_1 = &0, \; &b_2 =-&2, \; &b_3 = &2, \; &b_4 =-&2, \; &b_5 = &2, \; &b_6 = &0,  \\
        c_1 = &2, \; &c_2 = &2, \; &c_3 = &0, \; &c_4 = &0, \; &c_5 =-&2, \; &c_6 =-&2.  \\
    \end{array}
\end{align}
Similarly, $\mathcal{C}\odot \mathcal{K}$ satisfies
\begin{equation*}
    (\mathcal{C}\odot \mathcal{K})_{ijk} = (a_p i + b_p j + c_p k)h \lambda^{a_p i + b_p j + c_p k}, \quad \forall (i,j,k) \in D_{p}, \quad p=1,\cdots,6.
\end{equation*}

In this section, we first propose two fast implementations to accelerate the tensor-vector products $\mathcal{K} \times_i \bphi \times_j \bpsi$ in Subsection \ref{subsec:fast tensor vector product 1} and $(\mathcal{C} \odot \mathcal{K}) \times_i \bphi \times_j \bpsi$ in Subsection \ref{subsec:fast tensor vector product 2}, achieving a computation complexity of $O(N)$. Subsequently, we embed the fast tensor-vector product into the Sinkhorn algorithm, proposing the fast Sinkhorn algorithm with linear complexity for the MMOT problem in Subsection \ref{sec: FS3 algorithm}.

\subsection{Fast tensor-vector product of $\mathcal{K} \times_i \bphi \times_j \bpsi$}\label{subsec:fast tensor vector product 1}

The tensor-vector product writes
\begin{equation}\label{eqn:tensor-vector-product-1}
    (\mathcal{K} \times_i \bphi \times_j \bpsi)_k 
    = \sum\limits_{i=1}^N \sum\limits_{j=1}^N K_{ijk} \phi_i \psi_j
    = \sum_{p=1}^6 J_{k,p}, \quad k=1,2,\cdots,N,
\end{equation}    
in which $J_{k,p}$ are given by
\begin{align}\label{eqn: Jk}
    \begin{array}{l@{\quad}l}
    J_{k,1} = \sum\limits_{i=1}^k \sum\limits_{j=i}^k \phi_{i} \psi_{j} \lambda^{a_1 i + b_1 j +c_1 k}, & 
    J_{k,2} = \sum\limits_{j=1}^{k-1} \sum\limits_{i=j+1}^{k} \phi_{i} \psi_{j} \lambda^{a_2 i + b_2 j +c_2 k}, \\
    J_{k,3} = \sum\limits_{i=1}^k \sum\limits_{j=k+1}^{N} \phi_{i} \psi_{j} \lambda^{a_3 i + b_3 j +c_3 k}, & 
    J_{k,4} = \sum\limits_{j=1}^k \sum\limits_{i=k+1}^{N} \phi_{i} \psi_{j} \lambda^{a_4 i + b_4 j +c_4 k}, \\
    J_{k,5} = \sum\limits_{i=k+1}^N \sum\limits_{j=i}^{N} \phi_{i} \psi_{j} \lambda^{a_5 i + b_5 j +c_5 k}, & 
    J_{k,6} = \sum\limits_{j=k+1}^{N-1} \sum\limits_{i=j+1}^{N} \phi_{i} \psi_{j} \lambda^{a_6 i + b_6 j +c_6 k}.
    \end{array}
\end{align}
Directly computing Equation \eqref{eqn:tensor-vector-product-1} for all $k=1,\cdots,N$ takes $O(N^3)$ operations. In fact, $J_{k,p}$ satisfy the recurrence relations
\begin{align}\label{eqn:recurrence formula Jk}
    \begin{array}{l@{\quad}l}
        \left\{
        \begin{aligned}
            & J_{k,1} = \lambda^{c_1}J_{k-1,1} + \psi_{k} P_{k,1}, \\
            & J_{1,1} = \phi_1\psi_1 \lambda^{a_1 + b_1 + c_1}, \\
        \end{aligned}
        \right.
        & 
        \left\{
        \begin{aligned}
             & J_{k,2} = \lambda^{c_2}J_{k-1,2} + \phi_{k}P_{k-1,2}, \\
            & J_{1,2} = 0, \\
        \end{aligned}
        \right.   \\
        \\
        J_{k,3} =  P_{k,3}Q_{k+1,3},
        & J_{k,4} = P_{k,4}Q_{k+1,4}, \\
        \\
        \left\{
        \begin{aligned}
            & J_{k,5} = \lambda^{-c_5}J_{k+1,5} + \phi_{k+1} P_{k+1,5}, \\
            & J_{N,5} = 0, \\
        \end{aligned}
        \right.
        &
        \left\{
        \begin{aligned}
            & J_{k,6} = \lambda^{-c_6}J_{k+1,6} + \psi_{k+1} P_{k+2,6}, \\
            & J_{N,6} = 0,\;J_{N-1,6} = 0. \\
        \end{aligned}
        \right.
    \end{array}
\end{align}
Here $P_{k,p}$ and $Q_{k,p}$
have the forms of
\begin{align}\label{eqn: PkQk}
    \begin{array}{l@{\quad}l}
    P_{k,1} =\lambda^{(b_1 + c_1) k} \sum\limits_{i=1}^{k} \phi_i \lambda^{a_1 i}, &
    P_{k,2} = \lambda^{a_2(k+1) + c_2(k+1)} \sum\limits_{j=1}^{k} \psi_j \lambda^{b_2 j}, \\
    P_{k,3} = \lambda^{(b_3 + c_3) k} \sum\limits_{i=1}^k \phi_{i} \lambda^{a_3 i}, &
    Q_{k+1,3} = \lambda^{- b_3 k} \sum\limits_{j=k+1}^{N}  \psi_{j} \lambda^{b_3 j}, \\
    P_{k,4} = \lambda^{(a_4 + c_4) k} \sum\limits_{j=1}^k \psi_{j} \lambda^{b_4 j}, &
    Q_{k+1,4} = \lambda^{ - a_4 k} \sum\limits_{i=k+1}^{N}  \phi_{i} \lambda^{a_4 i}, \\
    P_{k,5} = \lambda^{a_5 k + c_5 (k-1)}\sum\limits_{j=k}^{N} \psi_j \lambda^{b_5 j}, &
    P_{k,6} =\lambda^{b_6 (k-1) + c_6 (k-2)}  \sum\limits_{i=k}^{N} \phi_i \lambda^{a_6 i}, 
    \end{array}
\end{align}
which can be also computed in a recursive manner:
\begin{align}\label{eqn:recurrence formula PkQk}
    \begin{array}{l@{\quad}l}
        \left\{
        \begin{aligned}
            & P_{k,1} = \lambda^{b_1 + c_1}P_{k-1,1} + \phi_k \lambda^{(a_1 + b_1 + c_1) k}, \\
            & P_{1,1} = \phi_1 \lambda^{a_1 + b_1 + c_1},
        \end{aligned}
        \right.
        & 
        \left\{
        \begin{aligned}
            & P_{k,2} = \lambda^{a_2 + c_2}P_{k-1,2} + \psi_k \lambda^{a_2(k+1) + b_2k + c_2(k+1)}, \\
            & P_{1,2} = \psi_1 \lambda^{2a_2 + b_2 + 2c_2},
        \end{aligned}
        \right.   \\
        \\
        \left\{
        \begin{aligned}
            & P_{k,3} = \lambda^{b_3 + c_3}P_{k-1,3} + \phi_k \lambda^{(a_3 + b_3 + c_3) k}, \\
            & P_{1,3} = \phi_1 \lambda^{a_3+b_3 + c_3}, \\
        \end{aligned}
        \right.
        &
        \left\{
        \begin{aligned}
            & Q_{k,3} = \lambda^{b_3} Q_{k+1,3} + \psi_k \lambda^{b_3} , \\
            & Q_{N,3} = \psi_N \lambda^{b_3},
        \end{aligned}
        \right. \\
        \\
        \left\{
        \begin{aligned}
            & P_{k,4} = \lambda^{a_4 + c_4}P_{k-1,4} + \psi_k \lambda^{(a_4 + b_4+c_4) k}, \\
            & P_{1,4} = \psi_1 \lambda^{a_4 + b_4+c_4}, \\
        \end{aligned}
        \right.
        &
        \left\{
        \begin{aligned}
            & Q_{k,4} = \lambda^{a_4} Q_{k+1,4} + \phi_k \lambda^{a_4}, \\
            & Q_{N,4} = \phi_N \lambda^{a_4},
        \end{aligned}
        \right. \\
        \\
        \left\{
        \begin{aligned}
            & P_{k,5} = \lambda^{-a_5-c_5}P_{k+1,5}+ \psi_k \lambda^{a_5k+b_5k+c_5 (k-1)}, \\
            & P_{N,5} = \psi_N \lambda^{a_5N+b_5N+c_5 (N-1)},
        \end{aligned}
        \right.
        &
        \left\{
        \begin{aligned}
            & P_{k,6} = \lambda^{-b_6-c_6}P_{k+1,6} + \phi_k \lambda^{a_6k+b_6 (k-1)+c_6 (k-2)}, \\
            & P_{N,6} = \phi_N \lambda^{a_6N+b_6 (N-1)+c_6 (N-2)}.
        \end{aligned}
        \right.
    \end{array}
\end{align}
The recurrence formulas \eqref{eqn:recurrence formula Jk} and \eqref{eqn:recurrence formula PkQk} induce the fast tensor-vector product algorithm (FTVP-1) to compute Equation \eqref{eqn:tensor-vector-product-1}. The pseudo-code is presented in Algorithm \ref{alg:fast tensor vector product}, with $O(N)$ computation complexity specified in Table \ref{tab:comutatioanl complexity of FTVP-1}.

\begin{algorithm}
	\caption{Fast tensor-vector product-1 (FTVP-1)}
	\label{alg:fast tensor vector product}   
	\hspace*{0.02in} {\bf Input:}  $\bphi$, $\bpsi$ of size $(N,1)$, $\lambda$ \\
	\hspace*{0.02in} {\bf Output:} 
	$\mathcal{K} \times_i \bphi \times_j \bpsi$
	\begin{algorithmic}[1] 
		\Procedure{FTVP-1}{$\bphi,\bpsi,\lambda$}  
		\State $\bJ_1,\bJ_2,\bJ_3,\bJ_4,\bJ_5,\bJ_6 = \boldsymbol{0}_{N}, J_{1,1}=\phi_1 \psi_1$
		\State $\boldsymbol{P}_{1},\boldsymbol{P}_{2},\boldsymbol{P}_{3},\boldsymbol{P}_{4},\boldsymbol{P}_{5},\boldsymbol{P}_{6},\boldsymbol{Q}_{3},\boldsymbol{Q}_{4} = \boldsymbol{0}_{N}$
		\State $P_{1,1}=\phi_1,P_{1,2}=\lambda^2\psi_1,P_{1,3}=\phi_1,Q_{N,3}=\lambda^{2}\psi_N$
		\State $P_{1,4}=\psi_1,Q_{N,4}=\lambda^{2}\phi_N,P_{N,5}=\lambda^2\psi_N,P_{N,6}=\lambda^4\phi_N$
            \For{$k = 1:N-1$}
		\State $P_{k+1,1} = \lambda^2 P_{k,1} + \phi_{k+1}$ 
		\State $P_{k+1,2} = \lambda^2 P_{k,2} + \lambda^2\psi_{k+1}$ 
            \State $P_{k+1,3} = \lambda^2 P_{k,3} + \phi_{k+1}$ 
            \State $Q_{N-k,3} = \lambda^2 Q_{N-k+1,3} + \lambda^2 \psi_{N-k}$ 
		\State $P_{k+1,4} = \lambda^2 P_{k,4} + \psi_{k+1}$
            \State $Q_{N-k,4} = \lambda^2 Q_{N-k+1,4} + \lambda^2 \phi_{N-k}$
		\State $P_{N-k,5} = \lambda^2 P_{N-k+1,5} + \lambda^2\psi_{N-k}$
		\State $P_{N-k,6} = \lambda^2 P_{N-k+1,6} + \lambda^4\phi_{N-k}$
		\EndFor
		\For{$k = 1:N-1$}
		\State $J_{k+1,1} = \lambda^2 J_{k,1} + \psi_{k+1} P_{k+1,1}$
		\State $J_{k+1,2} = \lambda^2 J_{k,2} + \phi_{k+1} P_{k,2}$
		\State $J_{k,3} = P_{k,3} Q_{k+1,3}$
		\State $J_{k,4} = P_{k,4} Q_{k+1,4}$
		\EndFor
		\For{$k=N:-1:2$}
		\State $J_{k-1,5} = \lambda^2 J_{k,5} +\phi_{k} P_{k,5}$
		\EndFor
		\For{$k = N-1:-1:2$}
		\State $J_{k-1,6} = \lambda^2 J_{k,6} + \psi_{k} P_{k+1,6}$
		\EndFor
		\\    \noindent\Return{$\bJ_1+\bJ_2+\bJ_3+\bJ_4+\bJ_5+\bJ_6$}
		\EndProcedure
	\end{algorithmic}  
\end{algorithm}  

\begin{table} [H]
	\centering
	\caption{The number of multiplicative and additive operations in FTVP-1.}
    \label{tab:comutatioanl complexity of FTVP-1}
	\begin{tabular}{c|ccc}
		\toprule
		 & $P_{k,p}$ & $Q_{k,p}$ & $J_{k,p}$  \\
        \midrule
        Number of operations & 18N &  6N & 14N \\
        Total & \multicolumn{3}{|c}{38N} \\
        \bottomrule
	\end{tabular}
\end{table}

\subsection{Fast tensor-vector product of $(\mathcal{C} \odot \mathcal{K}) \times_i \bphi \times_j \bpsi$
}\label{subsec:fast tensor vector product 2}

The tensor-vector product writes 
\begin{equation} \label{eqn:tensor-vector-product-2}
    \big((\mathcal{C} \odot \mathcal{K}) \times_i \bphi \times_j \bpsi\big)_k
    = \sum\limits_{i=1}^N \sum\limits_{j=1}^N c_{ijk}{K}_{ijk} \phi_i \psi_j
    = h \sum_{p=1}^6 \hat{J}_{k,p},
\end{equation}    
in which $\hat{J}_{k,p}$ for $p=1,2,\cdots,6$ are given by
\begin{align}\label{eqn: hat Jk}
    \begin{array}{l@{\quad}}
    \hat J_{k,1} = \sum\limits_{i=1}^k \sum\limits_{j=i}^k (a_1 i + b_1 j +c_1 k)\phi_{i} \psi_{j} \lambda^{a_1 i + b_1 j +c_1 k}, \\
    \hat J_{k,2} = \sum\limits_{j=1}^{k-1} \sum\limits_{i=j+1}^{k} (a_2 i + b_2 j + c_2 k) \phi_i \psi_j \lambda^{a_2 i + b_2 j + c_2 k}, \\
    \hat J_{k,3} = \sum\limits_{i=1}^{k} \sum\limits_{j=k+1}^{N} (a_3 i + b_3 j + c_3 k) \phi_i \psi_j \lambda^{a_3 i + b_3 j + c_3 k}, \\
    \hat J_{k,4} = \sum\limits_{j=1}^{k} \sum\limits_{i=k+1}^{N}  (a_4 i + b_4 j + c_4 k) \phi_i \psi_j \lambda^{a_4 i + b_4 j + c_4 k}, \\
    \hat J_{k,5} = \sum\limits_{i=k+1}^N \sum\limits_{j=i}^{N} (a_5 i + b_5 j + c_5 k)\phi_{i} \psi_{j} \lambda^{a_5 i + b_5 j +c_5 k}, \\ 
    \hat J_{k,6} = \sum\limits_{j=k+1}^N \sum\limits_{i=j+1}^{N} (a_6 i + b_6 j + c_6 k)\phi_{i} \psi_{j} \lambda^{a_6 i + b_6 j +c_6 k}.
    \end{array}
\end{align}
Directly computing Equation \eqref{eqn:tensor-vector-product-2} for all $k=1,\cdots,N$ takes $O(N^3)$ operations. In fact, $\hat J_{k,p}$ satisfy the recurrence relations
\begin{align}\label{eqn:recurrence formula hat Jk}
    \begin{array}{l}
        \left\{
        \begin{aligned}
            & \hat J_{k,1} = \lambda^{c_1}\big(\hat J_{k-1,1} + c_1 J_{k-1,1} \big) + \psi_k \big( R_{k,1} + (b_1 + c_1) k P_{k,1} \big), \\
            & \hat J_{1,1} = (a_1 + b_1 + c_1) \phi_1\psi_1 \lambda^{a_1 + b_1 + c_1},
        \end{aligned}
        \right.
        \\
        \\
        \left\{
        \begin{aligned}
            & \hat J_{k,2} = \lambda^{c_2} \big( \hat J_{k-1,2} + c_2 J_{k-1,2} \big) + \phi_k \big( R_{k-1,2} + (a_2 + c_2) k P_{k-1,2} \big), \\
            & \hat J_{1,2} = 0,
        \end{aligned}
        \right.
        \\
        \\
        \hat J_{k,3} = R_{k,3}Q_{k+1,3} +P_{k,3}S_{k+1,3} + c_3 k J_{k,3} ,
        \\
        \\
        \hat J_{k,4} = R_{k,4}Q_{k+1,4} + P_{k,4}S_{k+1,4} + c_4 k J_{k,4},
        \\
        \\
        \left\{
        \begin{aligned}
            & \hat J_{k,5} =  \lambda^{-c_5} \big( \hat J_{k+1,5} - c_5 J_{k+1,5} \big) +\phi_{k+1} \big( R_{k+1,5}+ (a_5(k+1) + c_5 k) P_{k+1,5} \big), \\
            & \hat J_{N,5} = 0,
        \end{aligned}
        \right.
        \\
        \\
        \left\{
        \begin{aligned}
            & \hat J_{k,6} =  \lambda^{-c_6} \big( \hat J_{k+1,6} - c_6 J_{k+1,6} \big) +\psi_{k+1} \big( R_{k+2,6}+ (b_6(k+1) + c_6k) P_{k+2,6} \big), \\
            & \hat J_{N,6} = 0, \quad \hat J_{N-1,6} = 0.
        \end{aligned}
        \right.
        \\
    \end{array}
\end{align}
Here $J_{k,p}, P_{k,p}, Q_{k,p}$ 
are the same as those in Equations \eqref{eqn: Jk} and \eqref{eqn: PkQk}. Additionally, $R_{k,p}$ and $S_{k,p}$
have the forms of
\begin{align}\label{eqn: RkSk}
    \begin{array}{l@{\quad}l}
    R_{k,1} = \lambda^{(b_1 + c_1)k}\sum_{i=1}^{k} a_1 i \phi_i \lambda^{a_1 i}, &
    R_{k,2} = \lambda^{(a_2+c_2)(k+1)}\sum\limits_{j=1}^k b_2 j \psi_j \lambda^{b_2 j}, \\
    R_{k,3} =  \lambda^{(b_3 + c_3) k}\sum\limits_{i=1}^{k} a_3 i \phi_i \lambda^{a_3 i}, &
    S_{k+1,3} =  \lambda^{- b_3 k}\sum\limits_{j=k+1}^{N} b_3 j \psi_j \lambda^{b_3 j}, \\
    R_{k,4} =  \lambda^{(a_4 + c_4) k}\sum\limits_{j=1}^{k} b_4 j \psi_j \lambda^{b_4 i}, &
    S_{k+1,4} =  \lambda^{-a_4 k} \sum\limits_{i=k+1}^{N} a_4 i \phi_i \lambda^{a_4 i}, \\
    R_{k,5} =  \lambda^{a_5 k + c_5 (k-1)}\sum\limits_{j=k}^N b_5 j \psi_j \lambda^{b_5 j}, &
    R_{k,6} = \lambda^{b_6(k-1)+c_6(k-2)}\sum\limits_{i=k}^N a_6 i \phi_i \lambda^{a_6 i}, 
    \end{array}
\end{align}
which can also be computed in a recursive manner:
\small{
\begin{align}\label{eqn:recurrence formula RkSk}
    \begin{array}{l@{\quad}l}
        \left\{
        \begin{aligned}
            & R_{k,1} = \lambda^{b_1 + c_1}R_{k-1,1} + a_1 k \phi_k \lambda^{(a_1+b_1+c_1) k}, \\
            & R_{1,1} = a_1 \phi_1 \lambda^{a_1+b_1+c_1},
        \end{aligned}
        \right.
        &
        \left\{
        \begin{aligned}
            & R_{k,2} = \lambda^{a_2 + c_2}R_{k-1,2} + b_2 k \psi_k \lambda^{a_2(k+1)+b_2k+c_2(k+1)}, \\
            & R_{1,2} = b_2\psi_1\lambda^{2a_2 + b_2 + 2c_2},
        \end{aligned}
        \right.
        \\
        \\
        \left\{
        \begin{aligned}
            & R_{k,3} = \lambda^{b_3 + c_3} R_{k-1,3} + a_3 k \phi_k \lambda^{(a_3+ b_3 + c_3) k}, \\
            & R_{1,3} = a_3 \phi_1 \lambda^{a_3+b_3+c_3},
        \end{aligned}
        \right.
        &
        \left\{
        \begin{aligned}
            & S_{k,3} = \lambda^{b_3} S_{k+1,3} + b_3 k \psi_k \lambda^{b_3}, \\
            & S_{N,3} = b_3 N \psi_N \lambda^{b_3},
        \end{aligned}
        \right.
        \\
        \\
        \left\{
        \begin{aligned}
            & R_{k,4} = \lambda^{a_4 + c_4}R_{k-1,4} + b_4 k \psi_k \lambda^{(a_4 + b_4+c_4) k}, \\
            & R_{1,4} = b_4 \psi_1 \lambda^{a_4 + b_4+c_4},
        \end{aligned}
        \right.
        &
        \left\{
        \begin{aligned}
            & S_{k,4} = \lambda^{a_4} S_{k+1,4} + a_4 k \phi_k \lambda^{a_4}, \\
            & S_{N,4} = a_4 N \phi_N \lambda^{a_4},
        \end{aligned}
        \right.
        \\
        \\
        \left\{
        \begin{aligned}
            & R_{k,5} = \lambda^{-a_5-c_5}R_{k+1,5} + b_5 k \psi_k \lambda^{a_5k+b_5 k+c_5(k-1)}, \\
            & R_{N,5} = b_5 \psi_N \lambda^{a_5N+b_5 N+c_5(N-1)},
        \end{aligned}
        \right.
        &
        \left\{
        \begin{aligned}
            & R_{k,6} = \lambda^{-b_6-c_6}R_{k+1,6} + a_6 k \phi_k \lambda^{a_6 k+b_6(k-1)+c_6(k-2)}, \\
            & R_{N,6} = a_6 \phi_N \lambda^{a_6 N+b_6(N-1)+c_6(N-2)}.
        \end{aligned}
        \right.
    \end{array}
\end{align}
}
The recurrence formulas \eqref{eqn:recurrence formula hat Jk} and \eqref{eqn:recurrence formula RkSk} induce the fast tensor-vector product algorithm (FTVP-2) to compute Equation \eqref{eqn:tensor-vector-product-2}. The pseudo-code is presented in Algorithm \ref{alg:fast tensor vector product-2}, with $O(N)$ computation complexity specified in Table \ref{tab:comutatioanl complexity of FTVP-2}.

\medskip

\begin{algorithm}
	\caption{Fast tensor-vector product-2 (FTVP-2)}
	\label{alg:fast tensor vector product-2}   
	\hspace*{0.02in} {\bf Input:}  $\bphi$, $\bpsi$ of size $(N,1)$, $\lambda,  h$ \\
	\hspace*{0.02in} {\bf Output:} 
    $(\mathcal{C} \odot \mathcal{K}) \times_i \bphi \times_j \bpsi$ 
	\begin{algorithmic}[1] 
		\Procedure{FTVP-2}{$\bphi,\bpsi,\lambda$}  
		\State $\hat{\bJ}_1,\hat{\bJ}_2,\hat{\bJ}_3,\hat{\bJ}_4,\hat{\bJ}_5,\hat{\bJ}_6 = \boldsymbol{0}_{N}$
		\State $\boldsymbol{R}_{1},\boldsymbol{R}_{2},\boldsymbol{R}_{3},\boldsymbol{R}_{4},\boldsymbol{R}_{5},\boldsymbol{R}_{6},\boldsymbol{S}_{3},\boldsymbol{S}_{4} = \boldsymbol{0}_{N}$
		\State $R_{1,1}=-2\phi_1,R_{1,2}=-2\lambda^2\psi_1,R_{1,3}=-2\phi_1,S_{N,3}=2N\lambda^{2}\psi_N$
		\State $R_{1,4}=-2\psi_1,S_{N,4}=2N\lambda^{2}\phi_N,R_{N,5}=2\lambda^2\psi_N,R_{N,6}=2\lambda^4\phi_N$
            \For{$k = 1:N-1$}
		\State $R_{k+1,1} = \lambda^2 R_{k,1} - 2(k+1)\phi_{k+1}$ 
		\State $R_{k+1,2} = \lambda^2 R_{k,2}-2(k+1)\lambda^2\psi_{k+1}$ 
            \State $R_{k+1,3} = \lambda^2 R_{k,3} - 2(k+1)\phi_{k+1}$ 
            \State $S_{N-k,3} = \lambda^2 S_{N-k+1,3} + 2(N-k) \lambda^2 \psi_{N-k}$ 
		\State $R_{k+1,4} = \lambda^2 R_{k,4} -2(k+1) \psi_{k+1}$
            \State $S_{N-k,4} = \lambda^2 S_{N-k+1,4} + 2(N-k) \lambda^2 \phi_{N-k}$
		\State $R_{N-k,5} = \lambda^2 R_{N-k+1,5} + 2(N-k)\lambda^2\psi_{N-k}$
		\State $R_{N-k,6} = \lambda^2 R_{N-k+1,6} + 2(N-k)\lambda^4\phi_{N-k}$
		\EndFor
		\For{$k = 1:N-1$}
		\State $\hat{J}_{k+1,1} = \lambda^2 \hat{J}_{k,1} + 2 \lambda^2 J_{k,1} +\psi_{k+1} R_{k+1,1} + 2(k+1)\psi_{k+1}P_{k+1,1}$
		\State $\hat{J}_{k+1,2} = \lambda^2 \hat{J}_{k,2} + 2 \lambda^2 J_{k,2} + \phi_{k+1} R_{k,2} + 2(k+1)\phi_{k+1}P_{k,2}$
		\State $\hat{J}_{k,3} = R_{k,3} Q_{k+1,3}+P_{k,3}S_{k+1,3}$
		\State $\hat{J}_{k,4} = R_{k,4} Q_{k+1,4}+P_{k,4}S_{k+1,4}$
		\EndFor
		\For{$k=N:-1:2$}
		\State $\hat{J}_{k-1,5} = \lambda^2 \hat{J}_{k,5} + 2 \lambda^2 J_{k,5}+\phi_{k} R_{k,5}-2(k-1)\phi_kP_{k,5}$
		\EndFor
		\For{$k = N-1:-1:2$}
		\State $\hat{J}_{k-1,6} = \lambda^2 \hat{J}_{k,6} + 2 \lambda^2 J_{k,6} + \psi_{k} R_{k+1,6}-2(k-1)\psi_{k} P_{k+1,6}$
		\EndFor
        \\
        \noindent\Return{$h(\hat{\bJ}_1+\hat{\bJ}_2+\hat{\bJ}_3+\hat{\bJ}_4+\hat{\bJ}_5+\hat{\bJ}_6)$}
		\EndProcedure
	\end{algorithmic}  
\end{algorithm}  

 \begin{table}
	\centering
	\caption{The number of multiplicative and additive operations in FTVP-2.}
    \label{tab:comutatioanl complexity of FTVP-2}
	\begin{tabular}{c|cccccc}
		\toprule
		 & $P_{k,p}$ & $Q_{k,p}$ & $J_{k,p}$ & $R_{k,p}$ & $S_{k,p}$ & $\hat J_{k,p}$ \\
        \midrule
        Number of operations & 18N &  6N & 14N & 30N & 10N & 48N \\
        Total & \multicolumn{6}{c}{126N} \\
        \bottomrule
	\end{tabular}
\end{table}

\subsection{Algorithm implementation} \label{sec: FS3 algorithm}

The major computational burden of the generalized Sinkhorn algorithm (Algorithm \ref{alg:Sinkhorn}) lies in the tensor-vector products related to $\mathcal{K}$ in lines 4-7 and that related to $\mathcal{C} \odot \mathcal{K}$ in line 8. Here we propose a step-by-step approach to reduce the computational complexity from $O(N^3)$ to $O(N)$:

\begin{itemize}
    \item $\mathcal{K} \times_j \bpsi \times_k \bvarphi$ in lines 4 and 7: we define $\mathcal{\tilde{K}} = (\tilde{K}_{ijk})$ satisfying $\tilde{K}_{ijk} = K_{kij} = \lambda^{b_p i + c_p j + a_p k},\;\forall (i,j,k) \in D_p$. Then it can be verified that 
    \begin{equation*}
        \mathcal{K} \times_j \bpsi \times_k \bvarphi = \mathcal{\tilde{K}} \times_i \bpsi \times_j \bvarphi,
    \end{equation*}
     which can be computed with $O(N)$ computational cost according to the discussion in Subsection \ref{subsec:fast tensor vector product 1}, because $\mathcal{\tilde{K}}$ has the similar form of Equation \eqref{eqn:form of K};

    \item $\mathcal{K} \times_i \bphi \times_k \bvarphi$ in lines 5 and 7: we can apply the same technique as used above to reduce the complexity to $O(N)$;

    \item $\mathcal{K} \times_i \bphi \times_j \bpsi$ in line 6: according to Algorithm \ref{alg:fast tensor vector product}, the complexity can be reduced to $O(N)$;

    \item $(\bC \odot \mathcal{K}) \times_i \bphi \times_j \bpsi$ in line 8: according to Algorithm \ref{alg:fast tensor vector product-2}, the complexity can be reduced to $O(N)$.
\end{itemize}
 
Based on the above discussions, we successfully accelerate the generalized Sinkhorn algorithm and propose a fast Sinkhorn algorithm with $O(N)$ computational complexity. The pseudo-code is presented in Algorithm \ref{alg:FS-3}.

\begin{algorithm}  
	\caption{The algorithm for MMOT with linear complexity}  
	\label{alg:FS-3}
	\hspace*{0.02in} {\bf Input:} $\bu, \bv, \bw$ of size $(N,1)$, $\varepsilon$, tol, itr\_max \\
	\hspace*{0.02in} {\bf Output:} $W_\epsilon(\bu,\bv,\bw)$
	\begin{algorithmic}[1] 
		\State initialize 
        $\bphi,\bpsi,\bvarphi = \frac{1}{N}\boldsymbol{1}_{N}$, 
        $t=0$, 
        $\mathrm{Res} = \mathrm{Inf}$
		\While{($t<$ itr\_max) \& ($\mathrm{Res} >$ tol)}
        \State $t \leftarrow t+1$
		\State $\bphi \leftarrow \bu \oslash \Call{FTVP-1}{\bpsi,\bvarphi,\lambda}$
		\State $\bpsi \leftarrow \bv \oslash \Call{FTVP-1}{\bphi,\bvarphi,\lambda}$
		\State $\bvarphi \leftarrow \bw \oslash \Call{FTVP-1}{\bphi,\bpsi,\lambda}$
		\State $\mathrm{Res} \leftarrow \mathrm{sum}\Big (\mid \bphi \odot  \Call{FTVP-1}{\bpsi,\bvarphi,\lambda} -\bu \mid +\mid \bpsi \odot \Call{FTVP-1}{\bphi,\bvarphi,\lambda}-\bv \mid \Big )$
		\EndWhile
          \State \Return{$W_\epsilon(\bu,\bv,\bw) = h \Big\langle \bvarphi,\Call{FTVP-2}{\bphi,\bpsi,\lambda} \Big\rangle.$}
	\end{algorithmic}  
\end{algorithm}

The log-domain stabilization can also be aggregated
into the fast tensor-vector product algorithm. The key point is to accelerate the tensor-vector products related to the rescaled tensor $\hat{\mathcal{K}} = (e^{\frac{\alpha_i + \beta_j + \gamma_k}{\varepsilon}} K_{ijk})$, which can be achieved by modifying the update formulas of $P_{k,p}, Q_{k,p}, J_{k,p}$ in FTVP-1. The pseudo-code is presented in Algorithm \ref{alg:fast tensor vector product with log}.

\begin{algorithm}
	\caption{Fast tensor-vector product with log-domain stabilization (FTVP-LOG)}
	\label{alg:fast tensor vector product with log}   
	\hspace*{0.02in} {\bf Input:} $\bphi$, $\bpsi$ of size $(N,1)$, $\lambda$ \\
	\hspace*{0.02in} {\bf Output:} 
	$\hat{\mathcal{K}} \times_i \bphi \times_j \bpsi$
	\begin{algorithmic}[1] 
		\Procedure{FTVP-1}{$\bphi,\bpsi,\lambda$}  
		\State $\bJ_1,\bJ_2,\bJ_3,\bJ_4,\bJ_5,\bJ_6 = \boldsymbol{0}_{N}, J_{1,1}=\phi_1 \psi_1$
		\State $\boldsymbol{P}_{1},\boldsymbol{P}_{2},\boldsymbol{P}_{3},\boldsymbol{P}_{4},\boldsymbol{P}_{5},\boldsymbol{P}_{6},\boldsymbol{Q}_{3},\boldsymbol{Q}_{4} = \boldsymbol{0}_{N}$
		\State $P_{1,1}=\phi_1,P_{1,2}=\lambda^2\psi_1,P_{1,3}=\phi_1,Q_{N,3}=\lambda^{2}\psi_N$
		\State $P_{1,4}=\psi_1,Q_{N,4}=\lambda^{2}\phi_N,P_{N,5}=\lambda^2\psi_N,P_{N,6}=\lambda^4\phi_N$
            \For{$k = 1:N-1$}
		\State $P_{k+1,1} = \lambda^2 {e^{\frac{\alpha_k-\alpha_{k+1}}{\varepsilon}}} P_{k,1} + \phi_{k+1}$ 
		\State $P_{k+1,2} = \lambda^2 {e^{\frac{\beta_k-\beta_{k+1}}{\varepsilon}}} P_{k,2} + \lambda^2\psi_{k+1}$ 
        \State $P_{k+1,3} = \lambda^2 {e^{\frac{\alpha_k-\alpha_{k+1}}{\varepsilon}}} P_{k,3} + \phi_{k+1}$ 
        \State $Q_{N-k,3} = 
    	\lambda^2 {e^{\frac{\beta_{N-k+1}-\beta_{N-k}}{\varepsilon}}} Q_{N-k+1,4} + \lambda^2 \psi_{N-k}$ 
		\State $P_{k+1,4} = \lambda^2 {e^{\frac{\beta_k-\beta_{k+1}}{\varepsilon}}} P_{k,2} + \psi_{k+1}$
        \State $Q_{N-k,4} = \lambda^2 {e^{\frac{\alpha_{N-k+1}-\alpha_{N-k}}{\varepsilon}}} Q_{N-k+1,4} + \lambda^2 \phi_{N-k}$
		\State $P_{N-k,5} = \lambda^2 {e^{\frac{\alpha_{N-k+1}-\alpha_{N-k}}{\varepsilon}}} P_{N-k+1,5} + \lambda^2 \psi_{N-k}$
		\State $P_{N-k,6} = \lambda^2 
    	{e^{\frac{\beta_{N-k+1}-\beta_{N-k}}{\varepsilon}}} P_{N-k+1,6} + \lambda^4 \phi_{N-k}$
		\EndFor
		\For{$k = 1:N-1$}
		\State $J_{k+1,1} = \lambda^2 {e^{\frac{\gamma_{k+1}-\gamma_{k}}{\varepsilon}}} J_{k,1} 
    	+ \psi_{k+1} {e^{\frac{\alpha_{k+1}+\beta_{k+1}+\gamma_{k+1}}{\varepsilon}}} P_{k+1,1}$
		\State $J_{k+1,2} = \lambda^2 {e^{\frac{\gamma_{k+1}-\gamma_{k}}{\varepsilon}}} J_{k,2} 
    	+ \phi_{k+1} {e^{\frac{\alpha_{k+1}+\beta_{k}+\gamma_{k+1}}{\varepsilon}}} P_{k,2}$
		\State $J_{k,3} = {e^{\frac{\alpha_k+\beta_{k+1}+\gamma_k}{\varepsilon}}} P_{k,3} Q_{k+1,3}$
		\State $J_{k,4} = {e^{\frac{\alpha_{k+1}+\beta_{k}+\gamma_k}{\varepsilon}}} P_{k,4} Q_{k+1,4}$
		\EndFor
		\For{$k=N:-1:2$}
		\State $J_{k-1,5} = \lambda^2 
    	{e^{\frac{\gamma_{k-1}-\gamma_{k}}{\varepsilon}}}
    	J_{k,5} + \phi_{k} {e^{\frac{\alpha_{k}+\beta_{k}+\gamma_{k-1}}{\varepsilon}}} P_{k,5}$
		\EndFor
		\For{$k = N-1:-1:2$}
		\State $J_{k-1,6} = \lambda^2 {e^{\frac{\gamma_{k-1}-\gamma_{k}}{\varepsilon}}} 
    	J_{k,6} + \psi_{k} {e^{\frac{\alpha_{k+1}+\beta_{k}+\gamma_{k-1}}{\varepsilon}}} P_{k+1,6}$
		\EndFor
		\\    \noindent\Return{$\bJ_1+\bJ_2+\bJ_3+\bJ_4+\bJ_5+\bJ_6$}
		\EndProcedure
	\end{algorithmic}  
\end{algorithm}

\section{Various Generalization}\label{sec:extension}

 In this section, we aim to extend the capabilities of our algorithm to higher dimensions and more margins.

\subsection{High dimension generalization}\label{sec:high dimension}

In this subsection, we consider the two-dimensional case of our algorithm as an example, noting that its generalization to three or higher dimensions follows a similar approach with no fundamental distinction. Similar to Section \ref{sec:mmot}, we consider the discrete two-dimensional MMOT problem, which solves the following minimization problem:
\begin{equation}\label{eqn:mot-2 problem}
	W(\bu,\bv,\bw)= \inf_{\mathcal{T} \in \Pi}\;\left \langle \mathcal{C}, \mathcal{T} \right \rangle =\inf_{\mathcal{T} \in \Pi}\;\sum\limits_{i_1=1}^{N_1} \sum\limits_{i_2=1}^{M_1} \sum\limits_{j_1=1}^{N_2} \sum\limits_{j_2=1}^{M_2} \sum\limits_{k_1=1}^{N_3} \sum\limits_{k_2=1}^{M_3} c_{i_1i_2j_1j_2k_1k_2} t_{i_1i_2j_1j_2k_1k_2}.
\end{equation}
Here $\mathcal{C}=(c_{i_1i_2j_1j_2k_1k_2}) \in \bbR^{N_1M_1\times N_2M_2 \times N_3M_3}$ is the cost tensor, and $\mathcal{T}=(t_{i_1i_2j_1j_2k_1k_2}) \in \bbR_{0+}^{N_1M_1\times N_2M_2 \times N_3M_3}$ is the multi-marginal transport plan, satisfying the linear constraints
\begin{multline}
    \label{eqn: transport map set-2}
	\mathcal{T} \in \Pi=\Big\{ \mathcal{T}=(t_{i_1i_2j_1j_2k_1k_2}) \Big|\;
        t_{i_1i_2j_1j_2k_1k_2} \ge 0,\;
	\sum\limits_{j_1=1}^{N_2} \sum\limits_{j_2=1}^{M_2} \sum\limits_{k_1=1}^{N_3} \sum\limits_{k_2=1}^{M_3} t_{i_1i_2j_1j_2k_1k_2}=u_{i_1i_2},\; \\
	\sum\limits_{i_1=1}^{N_1} \sum\limits_{i_2=1}^{M_1} \sum\limits_{k_1=1}^{N_3} \sum\limits_{k_2=1}^{M_3} t_{i_1i_2j_1j_2k_1k_2}=v_{j_1j_2},\; 
	\sum\limits_{i_1=1}^{N_1} \sum\limits_{i_2=1}^{M_1} \sum\limits_{j_1=1}^{N_2} \sum\limits_{j_2=1}^{M_2} t_{i_1i_2j_1j_2k_1k_2}=w_{k_1k_2}
	 \Big\}.   
\end{multline}

Similar to the discussion in Section \ref{sec:mmot}, our discussion is suitable for any $N_1$, $M_1$, $N_2$, $M_2$, $N_3$ and $M_3$. For the sake of simplicity, we assume $N_1 = N_2 = N_3 = N$ and $M_1 = M_2 = M_3 = M$ in the rest of the paper. Similar to the entropic regularization technique introduced in Section \ref{sec:mmot}, 
we consider the following regularized MMOT problem
\begin{equation*}
    W_\epsilon(\bu,\bv,\bw)
    =\inf_{\mathcal{T} \in \Pi}\;\sum\limits_{i_1,j_1,k_1=1}^{N} \sum\limits_{i_2,j_2,k_2=1}^{M} c_{i_1i_2j_1j_2k_1k_2} t_{i_1i_2j_1j_2k_1k_2} + \varepsilon t_{i_1i_2j_1j_2k_1k_2} \ln(t_{i_1i_2j_1j_2k_1k_2}).
\end{equation*}
Based on analogous derivation, the above problem achieves optimum when scaling variables $\bphi=(\phi_{i_1i_2}),\bpsi=(\psi_{j_1j_2}),\bvarphi =(\varphi_{k_1k_2}) \in \bbR^{NM}$ satisfy
\begin{align} 
	& \bphi \odot (\mathcal{K} \times_{j_1j_2} \bpsi  \times_{k_1k_2} \bvarphi) = \bu, \\  
	& \bpsi \odot (\mathcal{K} \times_{i_1i_2} \bphi \times_{k_1k_2} \bvarphi) = \bv, \\
	& \bvarphi \odot (\mathcal{K} \times_{i_1i_2} \bphi \times_{j_1j_2} \bpsi) = \bw, 
\end{align} 
in which the tensor-vector products write
\begin{align}
     (\mathcal{K} \times_{j_1j_2} \bpsi  \times_{k_1k_2} \bvarphi)_{i_1i_2} &= 
    \sum_{j_1,k_1=1}^N \sum_{j_2,k_2=1}^M K_{i_1i_2j_1j_2k_1k_2} \psi_{j_1j_2} \varphi_{k_1k_2}, \\
     (\mathcal{K} \times_{i_1i_2} \bphi \times_{k_1k_2} \bvarphi)_{j_1j_2} &= 
    \sum_{i_1,k_1=1}^N \sum_{i_2,k_2=1}^M K_{i_1i_2j_1j_2k_1k_2} \phi_{i_1i_2} \varphi_{k_1k_2}, \\
     (\mathcal{K} \times_{i_1i_2} \bphi \times_{j_1j_2} \bpsi)_{k_1k_2} &= \sum_{i_1,j_1=1}^N \sum_{i_2,j_2=1}^M K_{i_1i_2j_1j_2k_1k_2} \phi_{i_1i_2} \psi_{j_1j_2}.
\end{align}
The optimal transport plan is given by
\begin{equation*}
    t_{i_1i_2j_1j_2k_1k_2} = \phi_{i_1i_2}\psi_{j_1j_2}\varphi_{k_1k_2}K_{i_1i_2j_1j_2k_1k_2},
\end{equation*}
and the optimal transport distance is
\begin{multline}\label{eqn:mmot distance 2d}
    W_\epsilon(\bu,\bv,\bw) = \sum\limits_{i_1,j_1,k_1=1}^{N} \sum\limits_{i_2,j_2,k_2=1}^{M} c_{i_1 i_2j_1j_2k_1k_2} \phi_{i_1i_2}\psi_{j_1j_2} \varphi_{k_1k_2}  K_{i_1 i_2j_1j_2k_1k_2} \\
    = \Big\langle \bvarphi, (\mathcal{C} \odot \mathcal{K}) \times_{i_1i_2} \bphi \times_{j_1j_2} \bpsi \Big \rangle.
\end{multline}
Here, the cost tensor $C=(c_{i_1i_2j_1j_2k_1k_2}) \in \bbR^{NM \times NM \times NM}$ based on the $L^1$ norm of a uniform 2D mesh with a vertical spacing of $h_1$ and a horizontal spacing of $h_2$ writes
\begin{equation*}
    c_{i_1i_2j_1j_2k_1k_2}=(\abs{i_1-j_1}+\abs{i_1-k_1}+\abs{j_1-k_1})h_1+(\abs{i_2-j_2}+\abs{i_2-k_2}+\abs{j_2-k_2})h_2,
\end{equation*}
and the corresponding kernel tensor $\mathcal{K} = (K_{i_1i_2j_1j_2k_1k_2}) \in \bbR^{NM \times NM \times NM}$ writes
\begin{equation*}
K_{i_1i_2j_1j_2k_1k_2} = \lambda_1^{\abs{i_1-j_1}+\abs{i_1-k_1}+\abs{j_1-k_1}} \lambda_2^{\abs{i_2-j_2}+\abs{i_2-k_2}+\abs{j_2-k_2}},
\end{equation*}
    where $\lambda_1 = e^{-\frac{h_1}{\varepsilon}}$ and $\lambda_2 = e^{-\frac{h_2}{\varepsilon}}$.
The generalized Sinkhorn algorithm can be applied to calculate $(\bphi,\bpsi,\bvarphi)$ by iteratively updating
\begin{align}
	\bphi^{(t+1)} &=  \bu \oslash (\mathcal{K} \times_{j_1j_2} \bpsi^{(t)} \times_{k_1k_2} \bvarphi^{(t)}), 
    \label{eqn: tensor vector product 2D 1}\\
	\bpsi^{(t+1)} &=  \bv \oslash (\mathcal{K} \times_{i_1i_2} \bphi^{(t+1)} \times_{k_1k_2} \bvarphi^{(t)}),
    \label{eqn: tensor vector product 2D 2}\\
	\bvarphi^{(t+1)} &=  \bw \oslash (\mathcal{K} \times_{i_1i_2} \bphi^{(t+1)} \times_{j_1j_2} \bpsi^{(t+1)}). \label{eqn: tensor vector product 2D 3}
\end{align}

Similar to the discussion of our algorithm for the 1D problem in Section \ref{sec:FS-3}, the bottleneck of computational efficiency arises from the tensor-vector products, such as $\bJ = \mathcal{K} \times_{i_1i_2} \bphi \times_{j_1j_2} \bpsi$, whose element satisfies
\begin{equation}\label{eqn: Jk1k2}
     J_{k_1k_2} = \sum_{i_1,j_1=1}^N \sum_{i_2,j_2=1}^M \lambda_1^{\abs{i_1-j_1}+\abs{i_1-k_1}+\abs{j_1-k_1}} \lambda_2^{\abs{i_2-j_2}+\abs{i_2-k_2}+\abs{j_2-k_2}} \phi_{i_1i_2} \psi_{j_1j_2}.
\end{equation}
Obviously, directly calculating this tensor-vector product requires $O(N^3M^3)$ operations, which is an unacceptable computational cost. Now, we delve into a detailed discussion on reducing the computational cost to $O(NM)$.

First, $\bphi,\bpsi,\bJ \in \bbR^{NM}$ are flattened into 1D vectors in column-major order 
\begin{equation*}
     \bphi=(\bphi_{1},\bphi_{2},\cdots,\bphi_{M}), \quad \bpsi=(\bpsi_{1},\bpsi_{2},\cdots,\bpsi_{M}), \quad
     \bJ=(\bJ_{1},\bJ_{2},\cdots,\bJ_{M}),
\end{equation*}
in which the sub-vectors $\bphi_{i_2}, \bpsi_{j_2}, \bJ_{k_2} \in \bbR^{N}$ are
\begin{equation*}
\bphi_{i_2} = (\phi_{1i_2},\phi_{2i_2},\cdots,\phi_{Ni_2}), \quad
\bpsi_{j_2} = (\psi_{1j_2},\psi_{2j_2},\cdots,\psi_{Nj_2}), \quad
\bJ_{k_2} = (J_{1k_2},J_{2k_2},\cdots,J_{Nk_2}). 
\end{equation*}
By defining a tensor $\mathcal{K}_N = (K_{i_1j_1k_1,N}) \in \bbR^{N\times N\times N}$ with the element
\begin{equation*}
    K_{i_1j_1k_1,N} = \lambda_1^{\abs{i_1-j_1}+\abs{i_1-k_1}+\abs{j_1-k_1}},
\end{equation*}
the Equation \eqref{eqn: Jk1k2} can be rewritten in vector form as
\begin{equation*}
    \bJ_{k_2} = \sum_{i_2=1}^{M} \sum_{j_2=1}^{M}  \lambda_2^{\abs{i_2-j_2}+\abs{i_2-k_2}+\abs{j_2-k_2}} \big( \mathcal{K}_N \times_{i_1} \bphi_{i_2} \times_{j_1} \bpsi_{j_2} \big).
\end{equation*}
We use the same technique in Section \ref{sec:FS-3} to eliminate the absolute value operation above, separating the summations into 6 components according to the order of subscripts
\begin{equation*}
    \bJ_{k_2} = \sum_{p=1}^6 \bJ_{k_2,p},
\end{equation*}
each of which has a simple form of
\begin{align}\label{eqn: Jk form 2D}
	\begin{array}{l}
    	\bJ_{k_2,1} = \sum\limits_{i_2=1}^{k_2} \sum\limits_{j_2=i_2}^{k_2} \lambda_2^{2(k_2-i_2)} \big( \mathcal{K}_N \times_{i_1} \bphi_{i_2} \times_{j_1} \bpsi_{j_2} \big),
    	\\
    	\bJ_{k_2,2} = \sum\limits_{j_2=1}^{k_2-1} \sum\limits_{i_2=j_2+1}^{k_2} \lambda_2^{2(k_2-j_2)} \big( \mathcal{K}_N \times_{i_1} \bphi_{i_2} \times_{j_1} \bpsi_{j_2} \big),
    	\\
    	\bJ_{k_2,3} = \sum\limits_{i_2=1}^{k_2} \sum\limits_{j_2=k_2+1}^{M} \lambda_2^{2(j_2-i_2)} \big( \mathcal{K}_N \times_{i_1} \bphi_{i_2} \times_{j_1} \bpsi_{j_2} \big),
    	\\
    	\bJ_{k_2,4} = \sum\limits_{j_2=1}^{k_2} \sum\limits_{i_2=k_2+1}^{M} \lambda_2^{2(i_2-j_2)} \big( \mathcal{K}_N \times_{i_1} \bphi_{i_2} \times_{j_1} \bpsi_{j_2} \big),
    	\\
    	\bJ_{k_2,5} = \sum\limits_{i_2=k_2+1}^{M} \sum\limits_{j_2=i_2}^{M} \lambda_2^{2(j_2-k_2)} \big( \mathcal{K}_N \times_{i_1} \bphi_{i_2} \times_{j_1} \bpsi_{j_2} \big),
    	\\
    	\bJ_{k_2,6} = \sum\limits_{j_2=k_2+1}^{M-1} \sum\limits_{i_2=j_2+1}^{M} \lambda_2^{2(i_2-k_2)} \big( \mathcal{K}_N \times_{i_1} \bphi_{i_2} \times_{j_1} \bpsi_{j_2} \big).
    	\\
	\end{array}
\end{align}
In fact, $\bJ_{k_2,p}$ satisfy the recurrence relations
\begin{align*}
    \begin{array}{l@{\quad}l}
        \left\{
        \begin{aligned}
            & \bJ_{k_2,1} = \lambda_2^{2}\bJ_{k_2-1,1} + \mathcal{K}_N \times_{i_1} \bP_{k_2,1} \times_{j_1} \bpsi_{k_2}, \\
            & \bJ_{1,1} = \mathcal{K}_N \times_{i_1} \bphi_{1} \times_{j_1} \bpsi_{1},
        \end{aligned}
        \right.
        & 
        \left\{
        \begin{aligned}
             & \bJ_{k_2,2} = \lambda_2^{2}\bJ_{k_2-1,2} +  \mathcal{K}_N \times_{i_1} \bphi_{k_2} \times_{j_1} \bP_{k_2-1,2}, \\
            & \bJ_{1,2} = \boldsymbol{0},
        \end{aligned}
        \right.   \\
        \\
        \bJ_{k_2,3} =  \mathcal{K}_N \times_{i_1} \bP_{k_2,3} \times_{j_1} \boldsymbol{Q}_{k_2+1,3},
        &
        \bJ_{k_2,4} = \mathcal{K}_N \times_{i_1} \boldsymbol{Q}_{k_2+1,4} \times_{j_1} \bP_{k_2,4},
        \\
        \\
        \left\{
        \begin{aligned}
            & \bJ_{k_2-1,5} = \lambda_2^{2}\bJ_{k_2,5} + \mathcal{K}_N \times_{i_1} \bphi_{k_2} \times_{j_1} \bP_{k_2,5}, \\
            & \bJ_{N,5} = \boldsymbol{0}, 
        \end{aligned}
        \right.
        &
        \left\{
        \begin{aligned}
            & \bJ_{k_2-1,6} = \lambda_2^{2}\bJ_{k_2,6} + \mathcal{K}_N  \times_{i_1} \bP_{k_2+1,6} \times_{j_1} \bpsi_{k_2}, \\
            & \bJ_{N,6} = \boldsymbol{0},\;\bJ_{N-1,6} =\boldsymbol{0}. 
        \end{aligned}
        \right.
    \end{array}
\end{align*}
Here $\bP_{k_2,p}, \bQ_{k_2,p} \in \bbR^{N}$ have the forms of
\begin{align*}
    \begin{array}{l@{\quad}l}
    \bP_{k_2,1} = \sum\limits_{i_2 = 1}^{k_2} \lambda_2^{2(k_2-i_2)} \bphi_{i_2}, &
    \bP_{k_2,2} = \sum\limits_{j_2 = 1}^{k_2} \lambda_2^{2(k_2-j_2)} \bpsi_{j_2}, \\
    \bP_{k_2,3} = \sum\limits_{i_2 = 1}^{k_2} \lambda_2^{2(k_2-i_2)} \bphi_{i_2}, & 
    \bQ_{k_2,3} = \sum\limits_{j_2 = k_2}^{M} \lambda_2^{2(j_2-k_2)}\bpsi_{j_2},
    \\
    \bP_{k_2,4} = \sum\limits_{j_2 = 1}^{k_2} \lambda_2^{2(k_2-j_2)} \bpsi_{j_2}, & 
    \bQ_{k_2,4} = \sum\limits_{i_2 = k_2}^{M} \lambda_2^{2(i_2-k_2)}\bphi_{i_2}
    \\
    \bP_{k_2,5} = \sum\limits_{j_2 = k_2}^{M} \lambda_2^{2(j_2-k_2)}\bpsi_{j_2}, &
    \bP_{k_2,6} = \sum\limits_{i_2 = k_2}^{M} \lambda_2^{2(i_2-k_2)}\bphi_{i_2}, \\
    \end{array}
\end{align*}
which can also be computed in a recursive manner
\begin{align*}
    \begin{array}{l@{\quad}l}
        \left\{
        \begin{aligned}
            & \bP_{k_2,1} = \lambda_2^{2}\bP_{k_2-1,1} + \bphi_{k_2}, \\
            & \bP_{1,1} = \bphi_1,
        \end{aligned}
        \right.
        & 
        \left\{
        \begin{aligned}
            & \bP_{k_2,2} = \lambda_2^{2}\bP_{k_2-1,2} +  \lambda_2^{2}\bpsi_{k_2}, \\
            & \bP_{1,2} =  \lambda_2^{2}\bpsi_1,
        \end{aligned}
        \right.   \\
        \\
        \left\{
        \begin{aligned}
            & \bP_{k_2,3} = \lambda_2^{2} \bP_{k_2-1,3} + \bphi_{k_2}, \\
            & \bP_{1,3} = \bphi_1, \\
        \end{aligned}
        \right.
        &
        \left\{
        \begin{aligned}
            & \boldsymbol{Q}_{k_2,3} = \lambda_2^{2} \boldsymbol{Q}_{k_2+1,3} + \lambda_2^{2}\bpsi_{k_2} , \\
            & \boldsymbol{Q}_{N,3} =  \lambda_2^{2}\bpsi_N,
        \end{aligned}
        \right. \\
        \\
        \left\{
        \begin{aligned}
            & \bP_{k_2,4} = \lambda_2^{2} \bP_{k_2-1,4} + \bpsi_{k_2}, \\
            & \bP_{1,4} = \bpsi_1, \\
        \end{aligned}
        \right. 
        &
        \left\{
        \begin{aligned}
            & \boldsymbol{Q}_{k_2,4} = \lambda_2^{2} \boldsymbol{Q}_{k_2+1,4} + \lambda^{2}_2\bphi_{k_2}, \\
            & \boldsymbol{Q}_{N,4} = \lambda_2^{2 }\bphi_N,
        \end{aligned}
        \right. 
        \\
        \\
        \left\{
        \begin{aligned}
            & \bP_{k_2,5} = \lambda_2^{2}\bP_{k_2+1,5}+ \lambda_2^{2}\bpsi_{k_2} , \\
            & \bP_{N,5} = \lambda_2^{2}\bpsi_N,
        \end{aligned}
        \right.
        &
        \left\{
        \begin{aligned}
            & \bP_{k_2,6} = \lambda_2^{2}\bP_{k_2+1,6} + \lambda_2^{4}\bphi_{k_2} , \\
            & \bP_{N,6} = \lambda_2^{4}\bphi_N.
        \end{aligned}
        \right.
    \end{array}
\end{align*}
The pseudo-code of the implementation above is presented in Algorithm \ref{alg:fast tensor vector produc in two dimension}, with $O(MN)$ computational complexity. 
Similar idea can be naturally applied to accelerate the tensor-vector product $(\mathcal{C} \odot \mathcal{K}) \times_{i_1i_2} \bphi \times_{j_1j_2} \bpsi$ in Equation \eqref{eqn:mmot distance 2d}, and so are omitted.
\begin{algorithm}
	\caption{Fast tensor-vector product for 2D MMOT}
	\label{alg:fast tensor vector produc in two dimension}   
	\hspace*{0.02in} {\bf Input:}  $\bphi$, $\bpsi$ of size $(N \times M)$, $\lambda_1, \lambda_2$\\
	\hspace*{0.02in} {\bf Output:} 
	$\mathcal{K} \times_i \bphi \times_j \bpsi$
	\begin{algorithmic}[1] 
		\Procedure{FTVP\_{2d}}{$\bphi,\bpsi,\lambda_1,\lambda_2$}  
		\State $\bJ_{1,1}=\Call{FTVP-1}{\bphi_1,\bpsi_1,\lambda_1},\;\bJ_{1,2},\bJ_{1,3},\bJ_{1,4},\bJ_{1,5},\bJ_{1,6} = \boldsymbol{0}_{N}$
        \State$\bP_{1,1}=\bphi_1,\bP_{1,2}=\lambda_2^2\bpsi_1,\bP_{1,3}=\bphi_1,\boldsymbol{Q}_{M,3}=\lambda_2^{2}\bpsi_M$
		\State $\bP_{1,4}=\bpsi_1,\boldsymbol{Q}_{M,4}=\lambda_2^{2}\bphi_M,\bP_{M,5}=\lambda_2^2\bpsi_M,\bP_{M,6}=\lambda_2^4\bphi_M$
		\For{$k = 1:M-1$}
		\State $\bP_{k+1,1} = \lambda_2^2 \bP_{k,1} + \bphi_{k+1}$ 
		\State $\bP_{k+1,2} = \lambda_2^2 \bP_{k,2} + \lambda_2^2\bpsi_{k+1}$ 
            \State $\bP_{k+1,3} = \lambda_2^2 \bP_{k,3} + \bphi_{k+1}$ 
            \State $\boldsymbol{Q}_{M-k,3} = \lambda_2^2 \boldsymbol{Q}_{M-k+1,3} + \lambda_2^2 \bpsi_{M-k}$ 
		\State $\bP_{k+1,4} = \lambda_2^2 \bP_{k,4} + \bpsi_{k+1}$
            \State $\boldsymbol{Q}_{M-k,4} = \lambda_2^2 \boldsymbol{Q}_{M-k+1,4} + \lambda_2^2 \bphi_{M-k}$
		\State $\bP_{M-k,5} = \lambda_2^2 \bP_{M-k+1,5} + \lambda_2^2\bpsi_{M-k}$
		\State $\bP_{M-k,6} = \lambda_2^2 \bP_{M-k+1,6} + \lambda_2^4\bphi_{M-k}$
		\EndFor
		\For{$k = 1:M-1$}
		\State $\bJ_{k+1,1} = \lambda_2^2 \bJ_{k,1} + \Call{FTVP-1}{\bP_{k+1,1},\bpsi_{k+1},\lambda_1}$
		\State $\bJ_{k+1,2} = \lambda_2^2 \bJ_{k,2} + \Call{FTVP-1}{\bphi_{k+1},\bP_{k,2},\lambda_1}$
		\State $\bJ_{k,3} = \Call{FTVP-1}{\bP_{k,3}, \boldsymbol{Q}_{k+1,3},\lambda_1}$
		\State $\bJ_{k,4} = \Call{FTVP-1}{\boldsymbol{Q}_{k+1,4}, \bP_{k,4}, \lambda_1}$
		\EndFor
		\For{$k=M:-1:2$}
		\State $\bJ_{k-1,5} = \lambda_2^2 \bJ_{k,5} + \Call{FTVP-1}{\bphi_{k},\bP_{k,5},\lambda_1}$
		\EndFor
		\For{$k = M-1:-1:2$}
		\State $\bJ_{k-1,6} = \lambda_2^2 \bJ_{k,6} + \Call{FTVP-1}{\bP_{k+1,6},\bpsi_{k},\lambda_1}$
		\EndFor
        \For{$k = 1:M$}
        \State $\bJ_k = \bJ_{k,1} + \bJ_{k,2} + \bJ_{k,3} + \bJ_{k,4} + \bJ_{k,5} + \bJ_{k,6}$
        \EndFor
		\\
		\noindent\Return{$(\bJ_1, \bJ_2, \cdots, \bJ_M)$}
		\EndProcedure
	\end{algorithmic}  
\end{algorithm}  
\subsection{$l$-marginal generalization}\label{sec:k-marginal}

In this subsection, we extend our algorithm to solve general $l$-marginal optimal transport problems.
For ease of exposition, we focus on one-dimensional probabilistic distribution on a uniform mesh with a grid spacing of $h$.
Its generalization to high-dimensional distributions follows an analogous approach as discussed in Subsection \ref{sec:high dimension}.

Given $l$ discrete distributions $\bu_1,\bu_2,\cdots,\bu_l$ where $\bu_j=(u_{j,1},u_{j,2},\cdots,u_{j,N_j})\in \mbr^{N_j}$, the discrete $l$-marginal optimal transport problem \cite{benamou2015iterative} is formulated as 
\begin{equation}\label{eqn:k-mot problem}
	W(\bu_1,\bu_2,\cdots,\bu_l)= \inf_{\mathcal{T} \in \Pi}\;\left \langle \mathcal{C}, \mathcal{T} \right \rangle =\inf_{\mathcal{T} \in \Pi}\;\sum\limits_{i_1=1}^{N_1} \sum\limits_{i_2=1}^{N_2}\cdots \sum\limits_{i_l=1}^{N_l} c_{i_1 i_2\cdots i_l} t_{i_1 i_2\cdots i_l}.
\end{equation}
Here, $\mathcal{C}=(c_{i_1 i_2 \cdots i_l}) \in \bbR^{N_1 \times N_2\times \cdots \times N_l}$ is the $l$-th order cost tensor, and the $L_1$ norm based pairwise-interaction cost is
\begin{equation*}
    c_{i_1 i_2 \cdots i_l} = \sum_{1\le p<q \le l} \abs{i_p-i_q}h.
\end{equation*}
The collection of all admissible multi-marginal transport plan $\mathcal{T}=(t_{i_1 i_2 \cdots i_l}) \in \bbR_{0+}^{N_1 \times N_2\times \cdots \times N_l}$ is given by
\begin{equation}\label{eqn: k-marginal transport map set}
	 \Pi=\Big\{ \mathcal{T}=(t_{i_1 i_2 \cdots i_l}) \Big| t_{i_1 i_2 \cdots i_l} \ge 0, \; 
	\sum\limits_{i_1=1}^{N_1}\cdots \sum\limits_{i_{j-1}=1}^{N_{j-1}}  \sum\limits_{i_{j+1} =1}^{N_{j+1}}\cdots\sum\limits_{i_{l} =1}^{N_l} t_{i_1 i_2 \cdots i_l} =u_{j,i_j},\;\forall j=1,2,\cdots,l 
	 \Big\}.
\end{equation}

Without loss of generality, we assume $N_j = N$ for all $j=1,2,\cdots,l$.
By introducing the entropic regularization technique, the optimal transport plan $\mathcal{T}=(t_{i_1 i_2 \cdots i_l})$ of the regularized $l$-marginal optimal transport problem 
\begin{equation*}
	W_\varepsilon(\bu_1,\bu_2,\cdots,\bu_l)
 =\inf_{\mathcal{T} \in \Pi}\;\sum\limits_{i_1=1}^{N} \sum\limits_{i_2=1}^{N}\cdots \sum\limits_{i_l=1}^{N} c_{i_1 i_2\cdots i_l} t_{i_1 i_2\cdots i_l} + \varepsilon t_{i_1 i_2\cdots i_l} \ln{(t_{i_1 i_2\cdots i_l})},
\end{equation*}
writes
\begin{equation}\label{eqn: k-marginal T value}
	t_{i_1 i_2 \cdots i_l} =  \phi_{1,i_1}\phi_{2,i_2} \cdots \phi_{l,i_l} K_{i_1 i_2 \cdots i_l},\quad  \forall\;1\le i_1,i_2, \cdots, i_l \le N.
\end{equation}
Here, $\bphi_l = (\phi_{l,i_l}) \in \bbR_{0+}^N$ denotes the $l$-th scaling variables, and the corresponding kernel tensor $\mathcal{K} = (K_{i_1 i_2 \dots i_l}) \in \bbR^{N^l} $ is defined as
\begin{equation*}
    K_{i_1 i_2 \dots i_l} = \lambda^{\sum_{1\le p<q \le l} \abs{i_p-i_q}}, \quad \lambda = e^{-\frac{h}{\varepsilon}}. 
\end{equation*}

By substituting \eqref{eqn: k-marginal T value} into the marginal constraint \eqref{eqn: k-marginal transport map set}, the optimality condition \eqref{eqn: 3-marginal optimality} for the 3-marginal optimal transport problem can be analogously extended to general $l$-marginal cases, which is
\begin{equation}\label{eqn: k-marginal optimality}
	\bphi_j \odot \left(\mathcal{K} \times_{i_1}\bphi_1\cdots \times_{i_{j-1}}\bphi_{j-1}\times_{i_{j+1}}\bphi_{j+1} \cdots \times_{i_l}\bphi_{l} \right)  = \bu_j,\quad j=1,2,\cdots,l.
\end{equation}
The element of the above tensor-vector product has the form of
\begin{multline}
     (\mathcal{K} \times_{i_1}\bphi_1\cdots \times_{i_{j-1}}\bphi_{j-1}\times_{i_{j+1}}\bphi_{j+1} \cdots \times_{i_l}\bphi_{l})_{i_j} \\= 
    \sum_{i_1=1}^N \cdots \sum_{i_{j-1}=1}^N \sum_{i_{j+1}=1}^N\cdots \sum_{i_l=1}^N 
    K_{i_1 i_2 \cdots i_l} \phi_{1,i_1}\cdots \phi_{j-1,i_{j-1}}  \phi_{j+1,i_{j+1}}\cdots \phi_{l,i_l}.
\end{multline}
Hence, the corresponding Sinkhorn iteration scheme to find $\bphi_j$ writes
\begin{align}
	\bphi_1^{(t+1)} &=  \bu_1 \oslash \left(\mathcal{K} \times_{i_2}\bphi_2^{(t)} \times_{i_3}\bphi_3^{(t)} \cdots \times_{i_l}\bphi_l^{(t)} \right), \label{eqn:update formula 1 for k-mmot} \\
	\bphi_2^{(t+1)} &=  \bu_2 \oslash \left(\mathcal{K} \times_{i_1}\bphi_1^{(t+1)} \times_{i_3}\bphi_3^{(t)} \cdots \times_{i_l}\bphi_l^{(t)} \right), \label{eqn:update formula 2 for k-mmot} \\
                &\vdots \notag  \\
	\bphi_l^{(t+1)} &=  \bu_l \oslash \left(\mathcal{K} \times_{i_1}\bphi_1^{(t+1)} \times_{i_2}\bphi_2^{(t+1)} \cdots \times_{i_{l-1}}\bphi_{l-1}^{(t+1)} \right),\label{eqn:update formula 3 for k-mmot}
\end{align}
and the $l$-marginal optimal transport distance is given by
\begin{multline}\label{eqn:k-mmot distance}
    W_\epsilon(\bu_1,\bu_2,\cdots,\bu_l) = \sum_{i_1=1}^{N} \sum_{i_2=1}^{N}\cdots \sum_{i_l=1}^{N} c_{i_1 i_2\cdots i_l} \phi_{1,i_1}\phi_{2,i_2} \cdots \phi_{l,i_l}  K_{i_1 i_2\cdots i_l} \\
    = \Big\langle \bphi_l, (\mathcal{C} \odot \mathcal{K}) \times_{i_1} \bphi_1 \times_{i_2} \bphi_2 \cdots \times_{i_{l-1}} \bphi_{l-1} \Big \rangle.
\end{multline}

Similar to the discussion of $l=3$ in Section \ref{sec:FS-3}, the bottleneck of computational efficiency for the $l$-marginal optimal transport problem lies in the tensor-vector products
\begin{equation}\label{eqn:tensor-vector product 1 k-marginal}
    \bJ^{(l)} = \mathcal{K} \times_{i_1} \bphi_1 \times_{i_2} \bphi_2 \cdots \times_{i_{l-1}} \bphi_{l-1},
\end{equation}
and
\begin{equation}\label{eqn:tensor-vector product 2 k-marginal}
    \hat{\bJ}^{(l)} = (\mathcal{C} \odot \mathcal{K}) \times_{i_1} \bphi_1 \times_{i_2} \bphi_2 \cdots \times_{i_{l-1}} \bphi_{l-1},
\end{equation}
where the superscript $(l)$ indicates the $l$-marginals.
The elements of $\bJ^{(l)} = (J_{i_l}^{(l)})$ and $\hat{\bJ}^{(l)} = (\hat{J}_{i_l}^{(l)})$  have the forms of 
\begin{equation}\label{eqn: Jk k margin}
J_{i_l}^{(l)} = \sum\limits_{i_1=1}^N \cdots \sum\limits_{i_{l-1}=1}^N \lambda^{\sum_{1\le p<q \le l} \abs{i_p-i_q}} \prod_{q=1}^{l-1} \phi_{q,i_q},
\end{equation}
and
\begin{equation}\label{eqn: hat Jk k margin}
\hat{J}_{i_l}^{(l)} = h\sum\limits_{i_1=1}^N \cdots \sum\limits_{i_{l-1}=1}^N \Big(\sum_{1\le p<q \le l} \abs{i_p-i_q}\Big)\lambda^{\sum_{1\le p<q \le l} \abs{i_p-i_q}} \prod_{q=1}^{l-1} \phi_{q,i_q}.
\end{equation}
Obviously, the direct calculations require $O(N^{l})$ operations, which are unacceptable computational costs. 
We show as follows how to accelerate the tensor-vector products to linear computational complexity relative to support size $N$.

To start, we eliminate the absolute value operation in Equations \eqref{eqn: Jk k margin} and \eqref{eqn: hat Jk k margin} by separating the summations into $l!$ components according to the order of subscripts
\begin{equation*}
    \bJ^{(l)} = \sum\limits_{p=1}^{l!} \bJ_{p}^{(l)},
    \qquad
    \hat{\bJ}^{(l)} = h\sum\limits_{p=1}^{l!} \hat{\bJ}_{p}^{(l)},
\end{equation*}
in which the elements of $\bJ_{p}^{(l)} = (J_{i_l,p}^{(l)})$ and $\hat{\bJ}_{p}^{(l)} = (\hat{J}_{i_l,p}^{(l)})$ write
\begin{equation*}
    J_{i_l,p}^{(l)} = \sum\limits_{(i_1,\cdots,i_l) \in D_p}
    \Big( \lambda^{\alpha_{k,p} i_l} \prod_{q=1}^{l-1} \phi_{q,i_q} \lambda^{\alpha_{q,p} i_q}  \Big), 
\end{equation*}
and
\begin{equation*}
    \hat{J}_{i_l,p}^{(l)} = \sum\limits_{(i_1,\cdots,i_l) \in D_p}
    \Big( a_{l,p} i_l + \sum_{q=1}^{l-1} a_{q,p} i_q \Big)
    \Big( \lambda^{\alpha_{l,p} i_l} \prod_{q=1}^{l-1} \phi_{q,i_q} \lambda^{\alpha_{q,p} i_q}  \Big).
\end{equation*}
Here, $D_p$ are the sets defined by the order of subscripts
\begin{align*}
     D_1 &= \Big\{(i_1,\cdots,i_{l})\;\Big|\;1\le i_1 \le i_2 \le \cdots \le i_l \le N \Big\}, \\
     D_2 &= \Big\{(i_1,\cdots,i_{l})\;\Big|\;1\le i_2 < i_1 \le \cdots \le i_l \le N \Big\}, \\
    &\vdots \\
     D_{l!} &= \Big\{(i_1,\cdots,i_{l})\;\Big|\;1 \le i_l < i_{l-1} < \cdots < i_1 \le N \Big\}.
\end{align*}

We first consider the calculation of $\bJ_{p}^{(l)} = (J_{i_l,p}^{(l)})$ and $\hat{\bJ}_{p}^{(l)} = (\hat{J}_{i_l,p}^{(l)})$ for the case of $p=1$. 
It is straightforward to show that
\begin{equation}
    J_{i_l,1}^{(l)} 
    = \sum_{1 \le i_1 \le \cdots \le i_{l-1} \le i_{l}} 
    \Big( \lambda^{\alpha_{l,1} i_l} \prod_{q=1}^{l-1} \phi_{q,i_q} \lambda^{\alpha_{q,1} i_q}  \Big)
    = \lambda^{\alpha_{l,1}} J_{i_l-1,1}^{(l)} +  \lambda^{\alpha_{l,1} i_l} \phi_{l-1,i_l} J^{(l-1)}_{i_{l},1},
\end{equation}
in which
\begin{equation*}
    J^{(l-1)}_{i_{l},1} = \sum_{1 \le i_1 \le \cdots \le i_{l-2} \le i_l} \Big( \lambda^{\alpha_{l-1,1} i_l} \prod_{q=1}^{l-2} \phi_{q,i_q} \lambda^{\alpha_{q,1} i_q}  \Big).
\end{equation*}
The above formula suggests that to compute $\bJ_{1}^{(l)} = (J_{i_l,1}^{(l)})$ in the $l$-marginal case, only $3N$ multiplications and $N$ additions are required, provided that $\bJ_{1}^{(l-1)} = (J^{(l-1)}_{i_{l-1},1})$ has been obtained in the $(l-1)$-marginal case.
This recursive process with respect to $l$ can be iterated until reaching $\bJ_{1}^{(3)} = (J_{i_3,1}^{(3)})$ in the $3$-marginal case, enabling $\bJ_{1}^{(l)}$ to be computed with $O(N)$ operations.

Similarly, $\hat{J}_{i_l,1}^{(l)}$ satisfy
\begin{multline}
    \hat{J}_{i_l,1}^{(l)} 
    = \sum_{1 \le i_1 \le \cdots \le i_{l-1} \le i_{l}} 
    \Big( \alpha_{l,1} i_l + \sum_{q=1}^{l-1} \alpha_{q,1} i_q \Big)\Big( \lambda^{\alpha_{l,1} i_l} \prod_{q=1}^{l} \phi_{q,i_q} \lambda^{\alpha_{q,1} i_q}  \Big) \\
    = \lambda^{\alpha_{l,1}} \hat{J}_{i_l-1,1}^{(l)}
    + \lambda^{\alpha_{l,1}i_l} \phi_{l-1,i_l} \hat{J}_{i_l,1}^{(l-1)}
    + \alpha_{l,1} J_{i_l,1}^{(l)} 
    + \alpha_{l,1} \lambda^{\alpha_{l,1}i_l} \phi_{l-1,i_l} J_{i_l,1}^{(l-1)},
    \label{eqn:tensor-vector product cost p1}
\end{multline}
in which
\begin{equation*}
    \hat{J}_{i_l,1}^{(l-1)} = \sum_{1 \le i_1 \le \cdots \le i_{l-2} \le i_{l}} 
    \Big( \alpha_{l-1,1} i_l + \sum_{q=1}^{l-2} \alpha_{q,1} i_q \Big)\Big( \lambda^{\alpha_{l-1,1} i_l} \prod_{q=1}^{l-2} \phi_{q,i_q} \lambda^{\alpha_{q,1} i_q}  \Big).
\end{equation*}
Note that, in \eqref{eqn:tensor-vector product cost p1}, both of the computational complexity of $\bJ_1^{(l)}$ and $\bJ_1^{(l-1)}$ are $O(N)$. 
Therefore, to compute $\hat{\bJ}_{1}^{(l)} = (\hat{J}_{i_l,1}^{(l)})$ in the $l$-marginal case, only $9N$ multiplications and $3N$ additions are required if $\hat{\bJ}_{1}^{(l-1)} = (\hat{J}^{(l-1)}_{i_{l-1},1})$ in the $(l-1)$-marginal case has been obtained. 
This recursive process with respect to $l$ can be iterated until reaching $\hat{\bJ}_{1}^{(3)} = (\hat{J}_{i_3,1}^{(3)})$ in the $3$-marginal case,
thereby enabling the calculation of $\hat{\bJ}_{1}^{(l)}$ with $O(N)$ operations.

Similar recursive processes of $O(N)$ cost can be straightforwardly applied to the calculation of other $\bJ_{p}^{(l)}$ and $\hat{\bJ}_{p}^{(l)}$ and hence are not included here.
So far, the tensor-vector products \eqref{eqn:tensor-vector product 1 k-marginal} and \eqref{eqn:tensor-vector product 2 k-marginal} can both be calculated with linear computational complexity relative to support size $N$.

\section{Numerical  Experiments}\label{sec:numerical experiment}

In this section, we conduct numerical experiments to demonstrate the accuracy and efficiency of our algorithm in both 1D and 2D cases. The generalized Sinkhorn algorithm is regarded as the ground truth of the entropy regularized MMOT problem. In the following experiments, we illustrate the computational costs of the two algorithms, as well as the disparities in their results.

In the numerical experiments, the regularization parameter $\varepsilon$ is set to $10^{-1}$ for both our algorithm and the Sinkhorn algorithm unless otherwise specified. For the sake of efficiency comparison, we conduct 100 Sinkhorn iterations uniformly across all experiments. For each scenario, the computational time and relative error results are averaged over 100 experiments. All the experiments are conducted in Matlab R2020a on a PC with 16G RAM, 11th Gen Intel (R) Core (TM) i5-1135G7 CPU @2.40GHz.

\subsection{1D random distribution\label{sec:5.1}}

In the first experiment, we consider three discrete 1D random distributions on the interval $[0,1]$, sampled by $N$ grid points uniformly distributed. Three discrete vectors $\bu^1,\bu^2,\bu^3$ are randomly generated on the grid points, where each element follows a uniform distribution over $[0,1]$. We compare the performances of the original Sinkhorn algorithm and our algorithm in computing  $W_\varepsilon(\frac{u^1}{\Vert u^1 \Vert_1 },\frac{u^2}{\Vert u^2 \Vert_1 } ,\frac{u^3}{\Vert u^3 \Vert_1 })$.

\begin{table}[H]
	\centering
	\begin{tabular}{ccccl}
		\hline
		{\color[HTML]{333333} }                             & \multicolumn{2}{c}{\textbf{Computational time (s)}} &                                           & \multicolumn{1}{c}{}                                 \\ \cline{2-3}
		\multirow{-2}{*}{{\color[HTML]{333333} \textbf{N}}} & \textbf{Ours}          & \textbf{Sinkhorn}          & \multirow{-2}{*}{\textbf{Speed-up ratio}} & \multicolumn{1}{c}{\multirow{-2}{*}{$\Vert P_{O}-P_{S} \Vert_F$}} \\ \hline
		10                                                  & $4.99\times10^{-1}$               & $6.48\times10^0$                   & $1.30\times10^1$                                  & $1.13\times 10^{-16}$                                          \\
		20                                                  & $1.01\times10^0$               & $5.26\times10^1$                   & $5.18\times10^1$                                 & $5.53\times 10^{-17}$                                          \\
		40                                                  & $2.05\times10^0$               & $4.20\times10^2$                   & $2.05\times10^2$                                 & $5.68\times 10^{-17}$                                           \\
		80                                                  & $4.12\times10^0$              & $3.47\times10^3$                & $8.43\times10^2$                                 & $4.12\times 10^{-17}$                                          \\
		\multicolumn{1}{l}{160}                             & $8.27\times10^0$              & $2.74\times10^4$                          & $3.32\times10^3$                                         & \multicolumn{1}{c}{$3.58\times 10^{-17}$ }                                \\ \hline
	\end{tabular}
	\caption{The 1D random distribution problem. The comparison between the Sinkhorn algorithm
		and our algorithm with different numbers of grid points $N$. $P_{O}$ and $P_{S}$ are the transport plans of our algorithm and the Sinkhorn algorithm, respectively.}
	\label{table1}
\end{table}
	
\begin{figure}[htbp]
	\centering
	\includegraphics[width=0.8\textwidth]{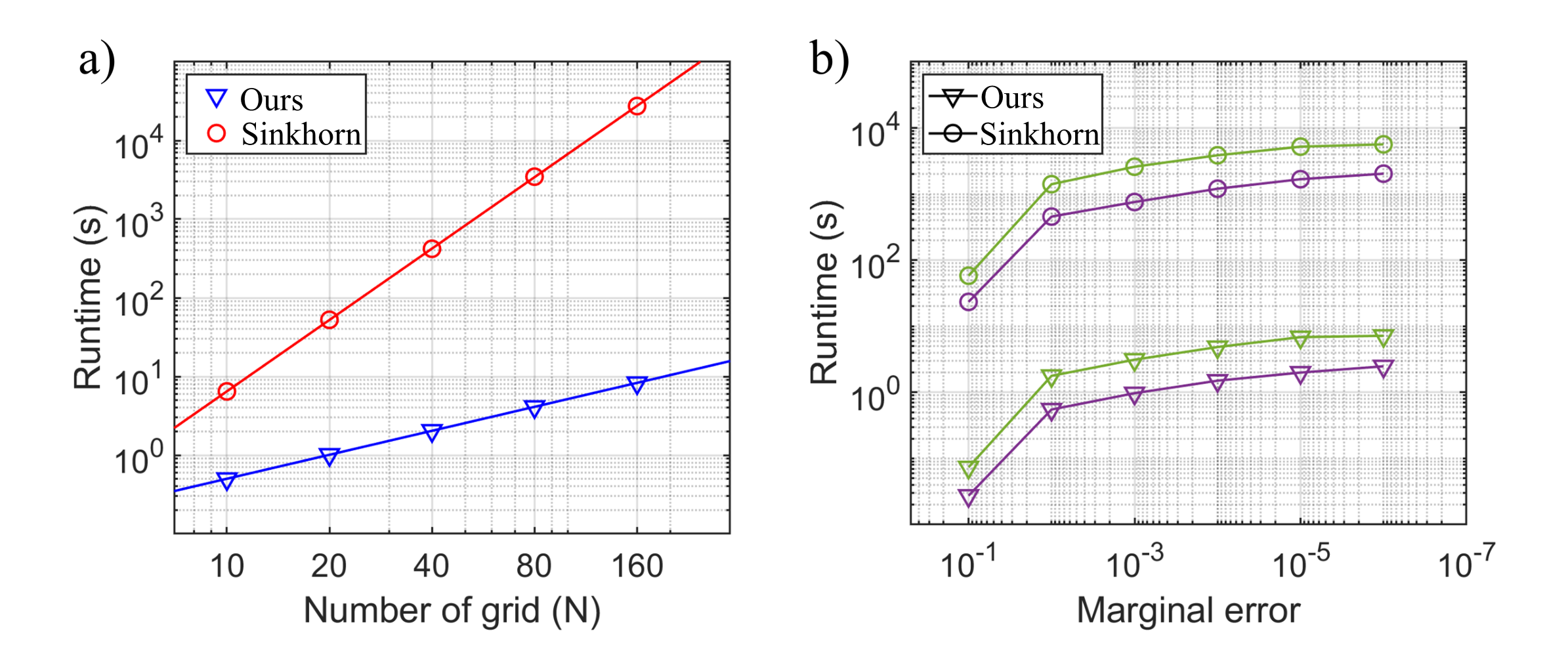}
	\caption{The 1D random distribution problem. (a): The computational time of our algorithm and the Sinkhorn algorithm with different numbers of grid points $N$. (b): The computational time required to reach the corresponding marginal error for $N = 80$ under $\varepsilon=0.1$ (purple) and $\varepsilon=0.05$ (green). }
	\label{fig1}
\end{figure}

In Table \ref{table1}, we output the computational time of the two algorithms and the difference between the transport plans computed by both algorithms. Notably, as the number of grid points increases, we observe a significant speed-up ratio while the optimal transport plans computed by both algorithms remain nearly identical. This shows the substantial advantage of our algorithm in computation efficiency. In Fig.~\ref{fig1}a, we compare the runtime of the Sinkhorn algorithm
and our algorithm with different numbers of grid points $N$. Employing linear regression, we obtain the empirical complexities of both algorithms: our algorithm exhibits $O(N^{1.01})$ complexity, and the Sinkhorn algorithm shows $O(N^{3.01})$ complexity. Additionally, Fig.~\ref{fig1}b illustrates the computational time required to achieve corresponding marginal errors for $N = 80$ under different $\varepsilon$, which emphasizes the overwhelming superiority of our algorithm in speed under different regularization parameters.

\subsection{The Ricker wavelet}\label{sec:5.2}

The Ricker wavelet is widely used in seismology for modeling the source time function, which has a form of
\begin{equation*}
	R(t)=A\left(1-2 \pi^{2} F_{0}^{2} t^{2}\right) e^{-\pi^{2} F_{0}^{2} t^{2}}.
\end{equation*}
Here $A$ denotes the wave amplitude, and $F_0$ is the dominant frequency. 

We evaluate the performances of the original Sinkhorn algorithm and our algorithm in measuring the discrepancy between three Ricker wavelets $f_1(t) = R(t-\tau_1)$, $f_2(t) = R(t-\tau_2), f_3(t) = R(t-\tau_3)$ based on the Wasserstein-1 distance. Without loss of generality, we set $A = 1, F_0 = 1$. To transform the Ricker wavelet to probability distribution function, we apply the normalization method introduced in \cite{li2022quadratic} and consider the following distance

\begin{equation}
	D(f_1, f_2, f_3)=W_{\varepsilon}\left(\frac{\frac{f_1^{2}}{\left\|f_1^{2}\right\|_1}+\delta}{1+L \delta}, \frac{\frac{f_2^{2}}{\left\|f_2^{2}\right\|_1}+\delta}{1+L \delta}, \frac{\frac{f_3^{2}}{\left\|f_3^{2}\right\|_1}+\delta}{1+L \delta}\right),
\end{equation}
where $L$ is used for normalization, and $\delta$ is a given parameter to improve numerical stability.

\begin{figure}[htbp]
	\centering
	\includegraphics[width=0.8\textwidth]{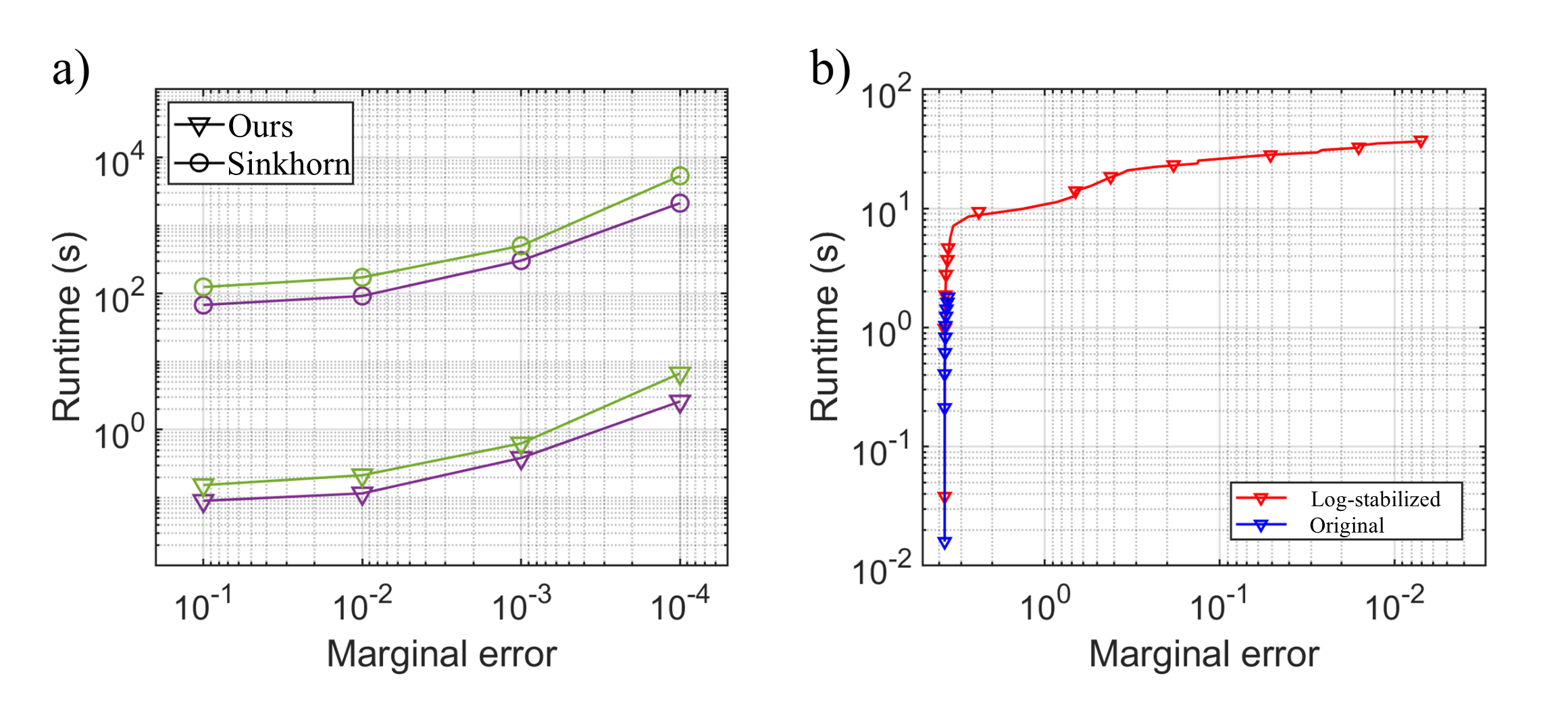}
	\caption{The Ricker wavelet problem.  (a): The computational time required to reach the corresponding marginal error for $N = 80$ under $\varepsilon=0.1$ (purple) and $\varepsilon=0.05$ (green). (b): The comparison between our algorithms with and without the log-domain stabilization for $N = 100, \varepsilon= 0.001$.}
	\label{fig2}
\end{figure}

In the experiment, we consider uniform grid points on the interval $[-2, 2]$, and set parameters $\tau_1 = 0,\;\tau_2 = 0.75,\;\tau_3 = 1.5,\;\delta = 10^{-3}$. In Table \ref{table2}, we output the computational time of both algorithms and the difference between the transport plans computed by both algorithms. The computational time required to achieve corresponding marginal errors for $N = 80$ under different $\varepsilon$ is shown in Fig.~\ref{fig2}a. As discussed in Subsection \ref{sec:5.1}, our algorithm maintains significant efficiency while preserving high computational precision compared to the Sinkhorn algorithm.

In Fig.~\ref{fig2}b, we show the importance of applying the log-domain stabilization technique when regularization parameters are relatively small, e.g., $\varepsilon=0.001$. Without the log-domain stabilization, our algorithm terminates at the 89th iteration due to exceeding the system threshold. However, with the log-domain stabilization applied, our algorithm maintains its efficiency without termination. 

\begin{table}[]
	\centering
	\begin{tabular}{ccccc}
		\hline
		{\color[HTML]{333333} }                             & \multicolumn{2}{c}{\textbf{Computational time (s)}} &                                           &                                  \\ \cline{2-3}
		\multirow{-2}{*}{{\color[HTML]{333333} \textbf{N}}} & \textbf{Ours}          & \textbf{Sinkhorn}          & \multirow{-2}{*}{\textbf{Speed-up ratio}} & \multirow{-2}{*}{{$\Vert P_{O}-P_{S} \Vert_F$}} \\ \hline
		10                                                  & $5.07\times10^{-1} $              & $6.78\times10^{0}  $                 & $1.34\times10^1 $                                 & $9.39\times10^{-17}$                \\
		20                                                  & $1.06\times10^0$               & $5.27\times10^1  $                 & $4.97\times10^{1}$                                  & $1.03\times10^{-16}  $                       \\
		40                                                  & $2.12\times10^0$               & $4.20\times10^2$                   & $2.68\times10^2$                                  & $1.93\times10^{-16}$                                \\
		80                                                  & $4.20\times10^0$               & $3.40\times10^3$                   & $8.10\times10^2$                                  & $7.36\times10^{-17}$                                \\
		\multicolumn{1}{l}{160}                             & $8.52\times10^0 $              & $2.71\times10^4$                   & $3.18\times10^3 $                                 & $1.16\times10^{-16}$                                \\ \hline
	\end{tabular}
	\caption{The Ricker wavelet problem. The comparison between the Sinkhorn algorithm and our algorithm with different numbers of grid points $N$. $P_{O}$ and $P_{S}$ are the transport plans of our algorithm and the Sinkhorn algorithm, respectively.}
	\label{table2}
\end{table}

\subsection{2D random distribution}\label{sec:5.3}

In this subsection, we evaluate the performance of the Sinkhorn algorithm and our algorithm in computing the Wasserstein-1 metric between three discrete 2D random distributions on $[0,1] \times [0,1]$. 
Those distributions are represented by $N \times M$ dimensional random vectors, where each element follows a uniform distribution over $[0,1]$. We use the basic settings mentioned in Subsection \ref{sec:high dimension}, and for simplicity, we set $M = N$.

\begin{figure}[htbp]
	\centering
	
	\includegraphics[width=0.8\textwidth]{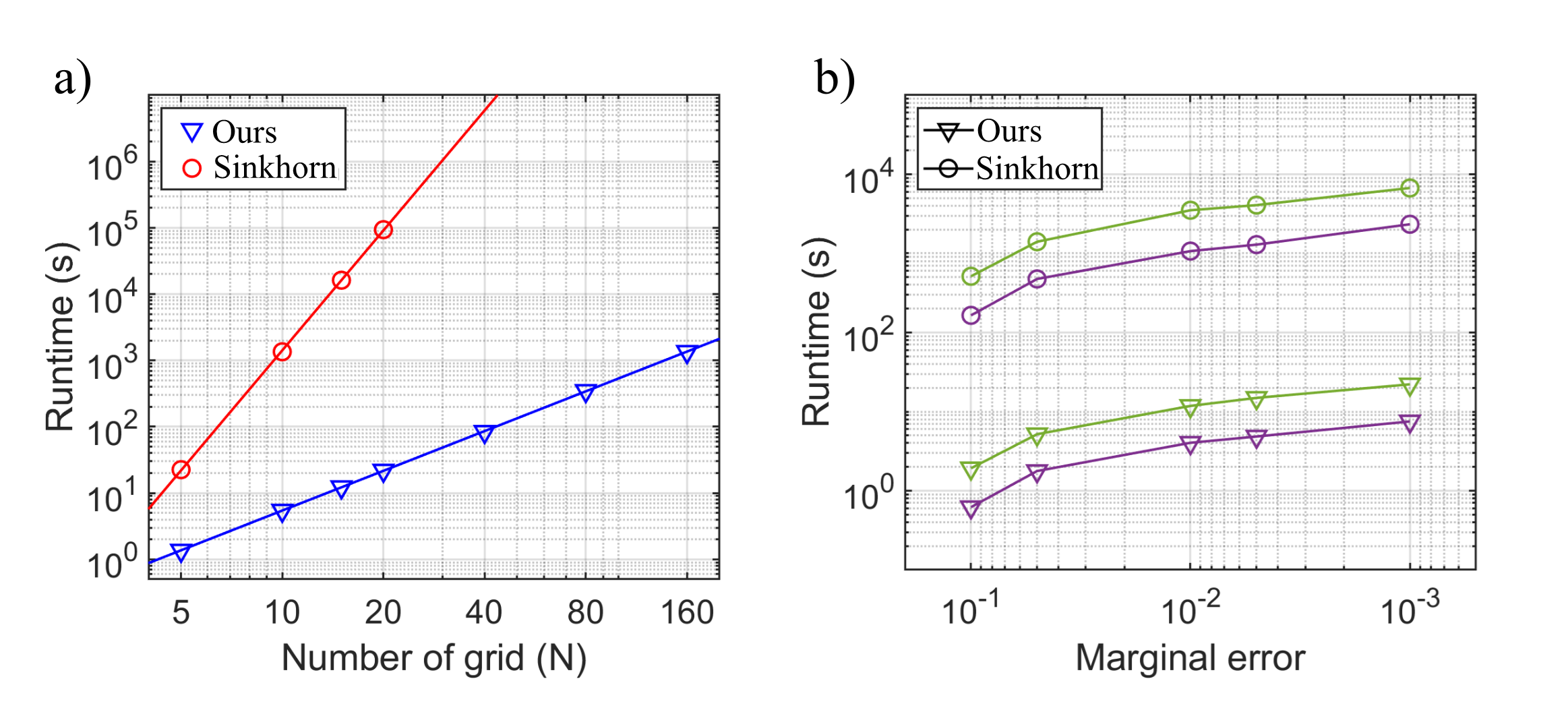}
	\caption{The 2D random distribution problem. (a): The comparison between the Sinkhorn algorithm
		and our algorithm with different numbers of grid points $N$. (b): The computational time required to reach the corresponding marginal error for $N = 10$ under $\varepsilon=0.1$ (purple) and $\varepsilon=0.01$ (green).}
	\label{fig3}
\end{figure}

Similar to Subsection \ref{sec:5.1}, we output the computational time of two algorithms and the difference between the transport plans computed by both algorithms in Table \ref{table3}. Additionally, Fig.~\ref{fig3}a illustrates the comparison between the Sinkhorn algorithm
and our algorithm with different numbers of grid points $N$. It is notable that with increasing $N$, the computation time of the Sinkhorn algorithm becomes unacceptable, which reaches several weeks. In contrast, the computational cost of our algorithm remains affordable, which draws the same conclusion as we discussed in the sections above. After conducting linear regression, we obtain that the empirical complexity of our algorithm is $O(N^{1.99})$ while that of the Sinkhorn algorithm is $O(N^{6.01})$. With the growth of dimension, the advantage of our algorithm in efficiency becomes more prominent. Furthermore, the computational time required to achieve corresponding marginal errors for $N = 10$ under different $\varepsilon$ is shown in Fig.~\ref{fig3}b. Without doubt, our algorithm outperforms the Sinkhorn algorithm by a wide margin in terms of speed under different regularization parameters.

\begin{table}[]
	\centering
	\begin{tabular}{ccccc}
		\hline
		{\color[HTML]{333333} }                               & \multicolumn{2}{c}{\textbf{Computational time (s)}} &                                           &                                  \\ \cline{2-3}
		\multirow{-2}{*}{{\color[HTML]{333333} \textbf{$N\times M$}}} & \textbf{Ours}          & \textbf{Sinkhorn}          & \multirow{-2}{*}{\textbf{Speed-up ratio}} & \multirow{-2}{*}{\textbf{$\Vert P_{O}-P_{S} \Vert_F$}} \\ \hline
		$10\times10$                                                 & $5.45 \times 10^0$               & $1.35\times10^3$                   & $2.48\times10^2$                                  & \multicolumn{1}{l}{$3.26\times 10^{-16}$}  \\
		$20\times20$                                                 & $2.19\times10^1$               & $9.36\times10^4$                          & $4.27\times10^3$                                         & $7.24\times10^{-16}$                                \\
		{$40\times40$}                             & $8.40\times10^1$               & -                          & -                                         & -                                \\ 
		{$80\times80$}                             & $3.48\times10^2$               & -                          & -                                         & -                                \\ 
		{$160\times160$}                             & $1.35\times10^3$               & -                          & -                                         & -\\ \hline
	\end{tabular}
	\caption{The 2D random distribution problem.The comparison between the Sinkhorn algorithm
		and our algorithm with different numbers of grid points $N$ and $\varepsilon = 0.1$. $P_{O}$ and $P_{S}$ are the transport plans of our algorithm and the Sinkhorn algorithm, respectively. We use `-' to denote the computational time exceeding $2\times10^6$s.}
	\label{table3}
\end{table}

\subsection{Multiple image matching problem }\label{sec:5.4}

\begin{figure}[htbp]
	\centering
	\includegraphics[scale=0.2]{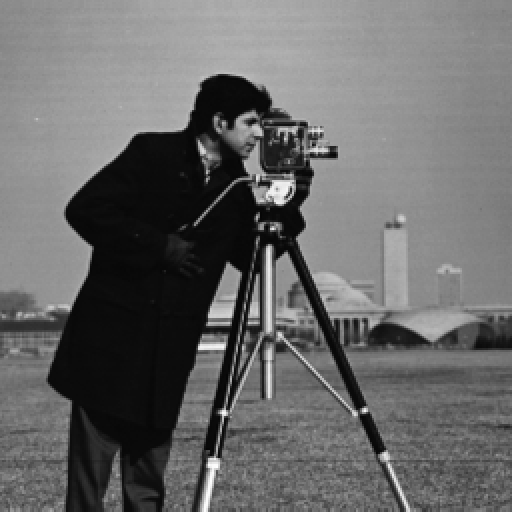}
	\includegraphics[scale=0.2]{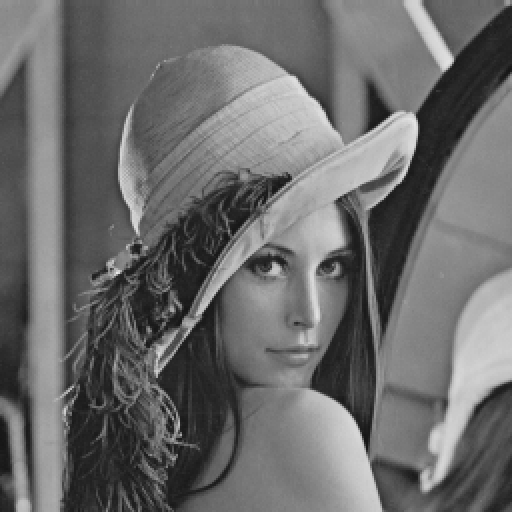}
	\includegraphics[scale=0.2]{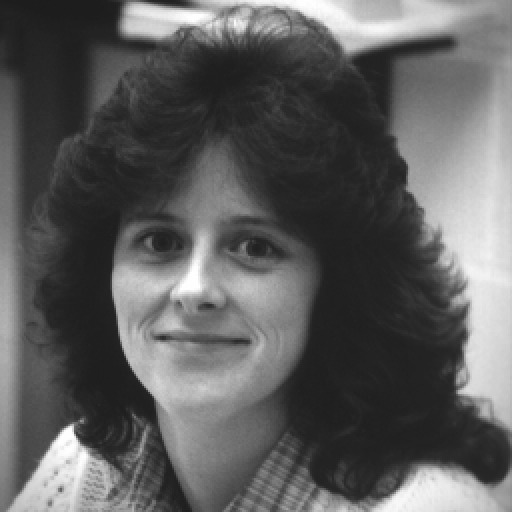}
	\caption{Illustration of images from DOTmark dataset.}
	\label{fig5}
\end{figure}

Nowadays, image matching problem has emerged as a trending topic in optimal transport \cite{burger2012regularized,maas2015generalized}. Our above experiments have already demonstrated the significant efficiency advantage of our algorithm over the traditional Sinkhorn algorithm, suggesting its considerable potential also in the realm of computer vision. Despite the extensive research on image matching, researchers are not satisfied with comparing two agents, but turn to multiple marginal matching problem which aims at learning transport plans to match a source domain to multiple target domains \cite{cao2019multi}. This multiple matching problem can be readily linked to the MMOT problem. Here we select three groups of images with different sizes from the DOTmark dataset which is specifically designed for discrete optimal transport \cite{DBLP:journals/access/SchrieberSG17} (see Fig.~\ref{fig5} for illustration). 

\begin{figure}[htbp]
	\centering
	\includegraphics[width=0.8\textwidth]{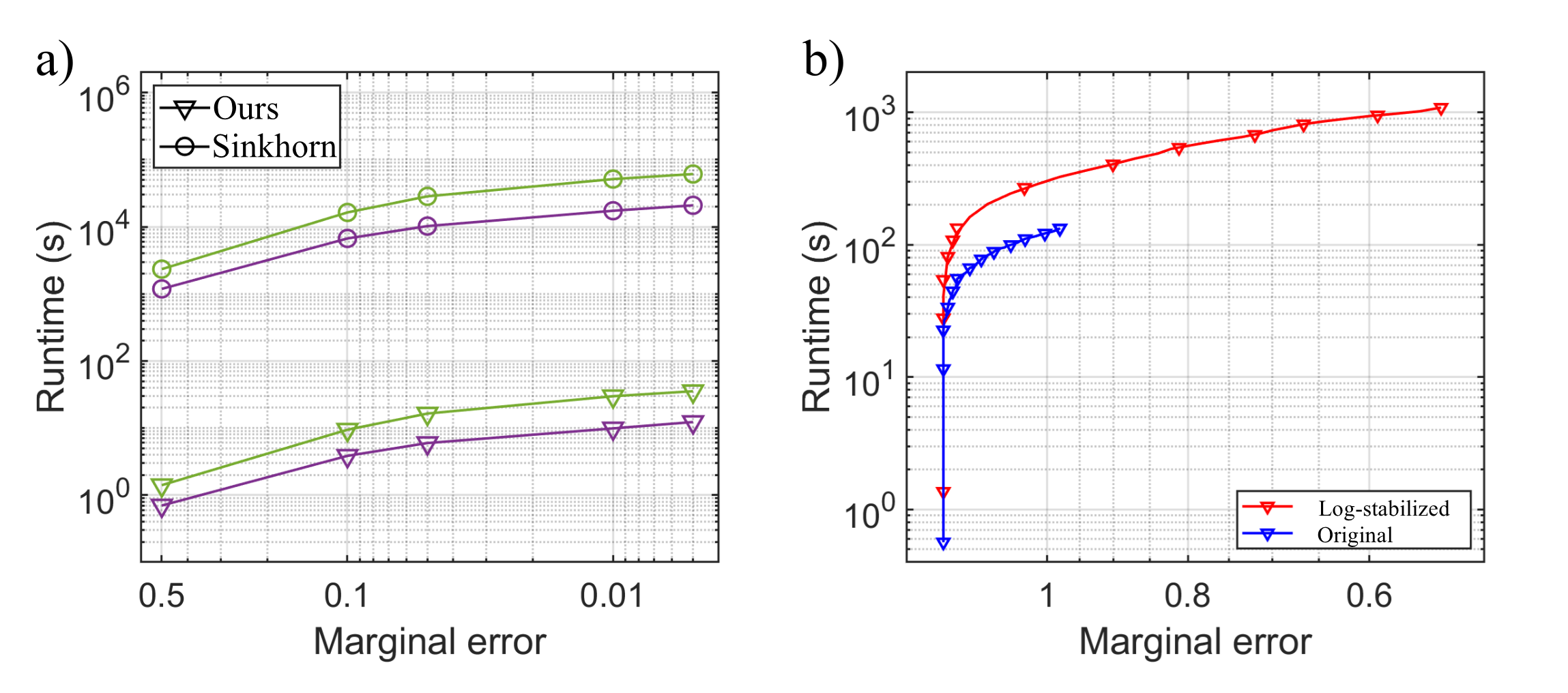}
	\caption{Multiple image matching problem.  (a): The computational time required to reach the corresponding marginal error for $N = 16$ under $\varepsilon=0.1$ (purple) and $\varepsilon=0.05$ (green). (b): The comparison between our algorithms with and without the log-domain stabilization for $N = 32, \varepsilon= 0.0005$.}
	\label{fig4}
\end{figure}

We transform images into 2D distributions on $[0,1] \times [0,1]$, represented by $M\times N$ dimensional vectors. Without loss of generality, we set $M = N$. The Wasserstain-1 metric between three images are computed by our algorithm and the original Sinkhron algorithm. 
In Table \ref{table4}, we present the computational time of the two algorithms and the difference between the transport plans computed by both algorithms. As mentioned in Subsection \ref{sec:5.3}, the computational time of the Sinkhorn algorithm soon becomes unbearable as $N$ increases. In contrast, the computational cost of our algorithm still remains affordable. The computational time required to achieve corresponding marginal errors for $N = 16$ under different $\varepsilon$ is shown in Fig.~\ref{fig4}a.

\begin{table}[H]
	\begin{tabular}{ccccc}
		\hline
		{\color[HTML]{333333} }                               & \multicolumn{2}{c}{\textbf{Computational time (s)}} &                                           &                                  \\ \cline{2-3}
		\multirow{-2}{*}{{\color[HTML]{333333} \textbf{$N\times M$}}} & \textbf{Ours}          & \textbf{Sinkhorn}          & \multirow{-2}{*}{\textbf{Speed-up ratio}} & \multirow{-2}{*}{\textbf{$\Vert P_{O}-P_{S} \Vert_F$}} \\ \hline
		$16\times16$                                               & $1.38\times10^1$               & $2.30\times10^5$                   & $1.67\times10^3 $                                 & \multicolumn{1}{l}{$6.41\times10^{-16}$}     \\
		$32\times32$                                               & $5.48\times10^1$               & $1.48\times10^6$                   & $2.70\times10^4 $                                 & \multicolumn{1}{l}{$1.06\times10^{-15}$}     \\
		$64\times64$                                               & $2.16\times10^2 $              & -                          & -                                         & -                                \\
		$128\times128 $                                            & $8.60\times10^2$               & -                          & -                                         & -                                \\
		$256\times256$                                             & $3.47\times10^3 $              & -                          & -                                         & -                                \\ \hline
	\end{tabular}
	\caption{Multiple image matching problem. The comparison between the Sinkhorn algorithm
		and our algorithm with numbers of grid points $N$ and $\varepsilon = 0.1$. $P_{O}$ and $P_{S}$ are the transport plans of our algorithm and the Sinkhorn algorithm, respectively. We use `-' to denote the computational time exceeding $2\times10^6$s.}
	\label{table4}
\end{table}

Finally, we show that the log-domain stabilization technique still works for two-dimensional problems with small regularization parameters. As shown in Fig.~\ref{fig4}b, our algorithm terminates at the 242nd iteration without this technique. In contrast, by employing the log-domain stabilization technique, our algorithm maintains its efficacy, proven to be an efficient method for computing the MMOT problem.

\section{Conclusion}\label{sec:conclusion}

In this paper, we propose an efficient numerical algorithm for solving the $L^1$ norm based multi-marginal optimal transport problem with linear complexity relative to support size $N$. This method accelerates the repeated tensor-vector products in the Sinkhorn algorithm by decomposing it into a summation of multiple components, each computed recursively with $O(N)$ additive and multiplicative operations. Furthermore, the log-domain stabilization technique is incorporated into our algorithm to avoid over- and underflow. In numerical experiments, we demonstrate that our algorithm achieves a significant speed advantage over the traditional method while maintaining accuracy. It reduces the computation cost for $l$-marginal case from $O(N^l)$ to $O(N)$ for any $l$. We anticipate that our algorithm will lead to a notable acceleration in various fields such as computer vision, machine learning, and transportation operations.

\section*{Acknowledgements}

This work was supported by National Natural Science Foundation of China Grant Nos. 12271289, 12031013 and 92270001, and Shanghai Municipal Science and Technology Major Project 2021SHZDZX0102. SJ was also supported by the
Fundamental Research Funds for the Central Universities.

\bibliographystyle{siam}
\bibliography{ref.bib}

\end{sloppypar}
\end{document}